\documentclass[reqno,11pt]{amsart}
\usepackage{amsmath, amsthm, a4, latexsym, amssymb}

\makeatletter\@addtoreset{equation}{}\makeatother

 \newfont{\bfit}{cmbxti10 scaled 1200}

 \newcommand{\eps}{\varepsilon}

 \newcommand{\R}{\mathbb{R}}
 
 \newcommand{\Z}{\mathbb{Z}}
 
 \newcommand{\prob}{\mathbb{P}}

 \newcommand{\me}{\mathbb{E}}
 
 \newcommand{\1}{{\sf 1}}

 \newcommand{\skrin}{{\mathcal N}}

 \newcommand{\heap}[2]{\genfrac{}{}{0pt}{}{#1}{#2}}
 \newcommand{\sfrac}[2]{\mbox{$\frac{#1}{#2}$}}

\newenvironment{Proof}[1]
{\vskip0.1cm\noindent{\bf #1}}{\vspace{0.15cm}}
\renewcommand{\subsection}{\secdef \subsct\sbsect}
\newcommand{\subsct}[2][default]{\refstepcounter{subsection}
\vspace{0.15cm}
{\flushleft\bf \arabic{section}.\arabic{subsection}~\bf #1  }
\nopagebreak\nopagebreak}
\newcommand{\sbsect}[1]{\vspace{0.1cm}\noindent
{\bf #1}\vspace{0.1cm}}

\newtheorem{theorem}{Theorem}

\newtheorem{lemma}[theorem]{Lemma}

\newtheorem{prop}[theorem]{Proposition}

\newcommand{\ssup}[1] {{\scriptscriptstyle{({#1}})}}
\newcommand{\ui}[1] {{\scriptscriptstyle{{#1}}}}

\theoremstyle{definition}
\newtheorem{remark}[theorem]{Remark}

\def\thebibliography#1{\section*{Bibliography}
  \list%
  {\arabic{enumi}.}
    {\settowidth\labelwidth{[#1]}\leftmargin\labelwidth
    \advance\leftmargin\labelsep
    \parsep0pt\itemsep0pt
    \usecounter{enumi}}
    \def\newblock{\hskip .11em plus .33em minus .07em}
    \sloppy                   
    \sfcode`\.=1000\relax}

\begin{document}
\title[Random walk in random scenery]
{Moderate deviations for random walk in random scenery}


\subjclass{Primary 60\thinspace F\thinspace10; Secondary
60\thinspace K\thinspace37}


\thanks{Corresponding author: Peter M\"orters}
\thanks{Running head: Random walk in random scenery}
 \maketitle \vspace{0.5cm}

\thispagestyle{empty}
\vspace{0.2cm}
\centerline{\textsc{Klaus Fleischmann, Peter M\"orters} and \textsc{Vitali Wachtel}}

\vspace{0.2cm}

\vspace{0.2cm}

\vfill
 \thispagestyle{empty}
 \setcounter{page}{0}

{\footnotesize

{\bf Klaus Fleischmann} and {\bf Vitali Wachtel}\\
Weierstrass Institute for Applied Analysis and Stochastics,\\
Mohrenstr.\ 39, D--10117 Berlin, Germany.\\
Email: \texttt{fleischm@wias-berlin.de}, resp. \texttt{vakhtel@wias-berlin.de}\\
URL: \texttt{http://www.wias-berlin.de/$\sim$fleischm}, resp. \texttt{http://www.wias-berlin.de/$\sim$vakhtel}\\

{\bf Peter M\"{o}rters}\\
University of Bath, Department of Mathematical Sciences,\\
Claverton Down, Bath BA2 7AY, United Kingdom.\\
Email: \texttt{maspm@bath.ac.uk}\\
URL: \texttt{http://people.bath.ac.uk/maspm/}\\}

\vfill

\newpage \vspace*{.1cm}
\begin{quote}{\small {\bf Abstract: } }
We investigate random walks in independent, identically distributed random sceneries
under the assumption that the scenery variables satisfy Cram\'er's  condition. We prove moderate
deviation principles in dimensions~$d\ge 2$, covering all those regimes where rate and
speed do not depend on the actual distribution of the scenery. In the case $d\ge 4$ we even obtain precise
asymptotics for the annealed probability of a moderate deviation, extending a classical central
limit theorem of Kesten and Spitzer. In $d\ge 3$, an important ingredient in the proofs are new
concentration inequalities for  self-intersection local times of random walks, which are of independent
interest, whilst in $d=2$ we use a recent moderate deviation result for  self-intersection local times,
which is due to Bass, Chen and Rosen.
\end{quote}
\vspace{0.5cm}



\section{Introduction}

In the world of stochastic processes in random environments, random walks in random scenery
represent a class of processes with fairly weak interaction. Nevertheless, they have deservedly
received a lot of attention since their introduction by Kesten and Spitzer~\cite{KS79} and,
independently, by Borodin~\cite{Bo79a, Bo79b}. A major reason for this interest is that in
$d\le 2$ the simple random walk in random scenery exhibits \emph{super-diffusive} behaviour.
However, in dimensions $d\ge 3$, when the underlying random walk visits most sites only once, the behaviour
of the random walk in random scenery is diffusive. Here finer features, like large deviation
behaviour, have to be studied in order to get an understanding of the interaction of walk and scenery.

To define random walk in random scenery, suppose
$\{S_n \colon n\ge 0\}$ is an underlying random walk on $\Z^d$ started at the origin, and
$\{\xi(z) \, : \, z \in\Z^d\}$ are independent, identically distributed
real-valued random variables, which are independent of the random walk
and which are called the scenery.  \emph{Random walk in random scenery} is the process $\{X_n \colon n\ge 0\}$ given by
$$X_n:= \sum_{1 \le k \le n} \xi(S_k)=\sum_{z\in\Z^d} \ell_n(z)\, \xi(z) \qquad \mbox{ for $n\ge 0$, }$$
where $\ell_n(z):=\sum_{1\le k\le n} \1\{S_k=z\}$ are the local times of the random walk
at the site~$z$.

Throughout this paper we make the following additional \emph{assumptions} on the model. The underlying walk is a
symmetric and aperiodic walk in dimensions $d\ge 2$, such that the covariance matrix~$\Gamma$ of $S_1$ is finite
and nondegenerate. Moreover, the random  variable $\xi(0)$ is centred, i.e. $\me \xi(0)=0$, with variance
$\sigma^2>0$, and satisfies $\mathbb{E}|\xi(0)|^3<\infty$ and \emph{Cram\'er's condition},
\begin{equation}\label{cram}
\me \big\{e^{\theta \xi(0) } \big\} <\infty \qquad \mbox{ for some } \theta>0. \\[2mm]
\end{equation}

The early papers by Kesten, Spitzer and Borodin establish \emph{central limit theorems} for
the random walk in random scenery. Indeed, it is (implicitly) shown in \cite{KS79} that, for $d\ge 3$,
\begin{equation}\label{KS}
\frac{X_n}{\sqrt{n}} \,\stackrel{n\uparrow\infty}{\Longrightarrow} \, \skrin\big(0, \sigma(2G(0)-1)\big),
\end{equation}
where $G$ is the Green's function of the underlying random walk.
Bolthausen in \cite{Bo89} extended this to the planar case by showing that
$$\frac{X_n}{\sqrt{n \log n}}\,\stackrel{n\uparrow\infty}{\Longrightarrow} \, \skrin(0, \pi^{-1}).$$
Hence, \emph{moderate and large deviation problems} for the random walk in random scenery
deal with the asymptotic behaviour of $\prob\{X_n \ge b_n\} $ for $b_n \gg \sqrt{n}$,
i.e. $\lim b_n/\sqrt{n}=\infty$, if $d\ge 3$, and $b_n \gg \sqrt{n \log n}$ if $d=2$.
Let us remark for completeness that Kesten and Spitzer
have also established a limit theorem in distribution for $X_n/n^{3/4}$ with non-Gaussian limits for $d=1$,
a case we do not consider in this paper as large and moderate deviations are more or less fully understood
in this case.\footnote{This information was communicated to us by F.~Castell.}

Large deviation problems for random walks in random scenery in
dimensions $d\ge 2$ have only recently attracted attention,
see~\cite{GP02, GHK06, GKS05, As06, AC05, AC06}, and also
\cite{CP01,AC03, Ca04} where Brownian motions are used in place of
random walks. The fascination of this subject stems {f}rom the rich
behaviour that comes to light when large deviations are
investigated. The intricate interplay of the walk with the scenery
leads to a large number of different regimes depending on
\begin{itemize}
\item the dimension $d$ of the underlying lattice $\Z^d$,
\item the upper tail behaviour of the scenery variable,
\item the size of the deviation studied,
\end{itemize}
to name just the most important ones. For example, Asselah and Castell~\cite{AC06},
restricting attention to dimensions~$d\ge 5$ and scenery variables with superexponential
decay of upper tails, have identified \emph{five} regimes with different large deviation speeds.
Heuristically, in each regime the walk and the scenery `cooperate' in a different way
to obtain the deviating behaviour. Up to now only one of these regimes has been fully
treated, including the discussion of explicit rate functions. This is
the \emph{very large deviation} regime discussed (together with a number of boundary cases)
by Gantert, K\"onig and Shi in~\cite{GKS05}. In this regime it is assumed that
$$\log \prob\{ \xi(0) > x\} \sim -D\, x^q \qquad\mbox{ as } x\uparrow\infty,$$
for some $D>0$ and $q>d/2$. Then, for any $n\ll b_n \ll n^{\frac{1+q}q}$, as $n\uparrow\infty$,
\begin{equation}\label{NinaWolf}
\log\prob\big\{ X_n > b_n \big\} \sim K\, n^{-\frac{2q-d}{d+2}}\,b_n^{\frac{2q}{d+2}},
\end{equation}
where $K=K(D,q,d)>0$ is a constant given explicitly in terms of a variational problem. The underlying
strategy is that the random walk \emph{contracts} to grow at a speed of $$n^{\frac{1+q}{d+2}}/b_n^{\frac{q}{d+2}}
\ll n^{\frac 12},$$ and the scenery adopts values of
size $b_n/n$ on the range of the walk.
The right hand side in \eqref{NinaWolf} represents the combined cost of these
two deviations.
\smallskip

In the present paper we study \emph{moderate deviation principles},
providing a full analysis including explicit rate functions and, in
dimensions $d\ge 4$, even exact asymptotics of moderate deviation
probabilities. We consider as moderate deviations the regimes
extending {f}rom the central limit scaling up to the point where
either the deviation speed or the rate function start to depend on
the actual distribution of the scenery, or in other words where tail
conditions stronger than Cram\'er's condition would have an impact
on the speed or rate of the deviations.

Heuristically, our results, which will be described in detail in the next section, show that in $d\ge 3$
throughout the moderate deviation regime the deviation is achieved by a moderate deviation of the scenery
without any contribution {f}rom the walk. The rates therefore agree with those obtained for \emph{fixed} walk
in a random scenery by Guillotin-Plantard in~\cite{GP02}. Crucial ingredients of our proofs are concentration
inequalities for self-intersection  local times of random walks, see Proposition~\ref{centredmoment}.
Our exact asymptotic results for the moderate deviation probabilities build on classical ideas of Cram\'er.
\medskip

In $d=2$, by contrast, the moderate deviation regime splits in two parts. If  $\sqrt{n \log n} \ll
b_n\ll \sqrt{n}\log n$ then, again,  we only have a contribution {f}rom the scenery and
the walk exhibits typical behaviour. However, if $\sqrt{n} \log n\ll b_n \ll  n/\log n$
the random walk \emph{contracts}, though in a much more delicate way than in the very large
deviation regime: The self-intersection local times
of the walk, which normally are of order $n \log n$ are now increased to be of
order $\sqrt{n} b_n$. At the same time, on the (contracted) range of the walk,
the scenery values perform a moderate deviation and take values of size $b_n/n$.
Our results in the case $d=2$ rely on moderate deviation principles for renormalised
self-intersection local times of planar random walks recently obtained by Bass,
Chen and Rosen~\cite{BCR06}.

\section{Main results}

Recall that we assume that the random variable $\xi(0)$ satisfies Cram\'er's condition~\eqref{cram}
and $\sigma^2>0$ denotes its variance. For $d\ge 3$ we define the Green's function of the random walk by
$$G(x):=\sum_{k=0}^\infty \prob\{ S_k=x\} \qquad\mbox{ for $x\in\Z^d$.}$$

\begin{theorem}[Refined moderate deviations in dimensions $d\ge 4$]\label{maind>3}\ \\ \vspace{-.5cm}

There exists a regularly varying sequence $(a_n)$ of index $\frac 23$, such that,
if $d\ge 4$ and $n^{\frac 12} \ll b_n \ll a_n$, then
$$\prob\big\{ X_n \ge b_n \big\} \sim 1 - \Phi\bigg(\frac{b_n}{\sqrt{\sigma^2\,n\,(2G(0)-1)}}\bigg)
\qquad\mbox{ as $n\uparrow\infty$}  \, ,$$
where $\Phi$ denotes the standard normal distribution function.
\end{theorem}
\smallskip

\begin{remark}
This result extends the central limit theorem~\eqref{KS} to the
moderate deviation regime. Note that asymptotics of this degree of precision are very rarely
encountered in stochastic processes beyond the independent case.
In this theorem we are restricted to dimensions $d\ge 4$
as our proof requires an analysis of \emph{triple} self-intersections of random walks,
for which $d=3$ is the critical dimension.
\end{remark}

In dimension~$d=3$ we can no longer provide \emph{precise} asymptotics, but we can still prove a full
moderate deviation principle with the same speed and rate function as in $d\ge 4$.

\begin{theorem}[Moderate deviations in dimensions $d\ge 3$]\label{maind>2}\ \\ \vspace{-.5cm}

If $d\ge 3$ and $n^{\frac 12} \ll b_n \ll n^{\frac 23}$, then, as $n\uparrow\infty$,
$$\log \prob\big\{ X_n \ge b_n \big\} \sim - \frac{b_n^2}{n} \, \frac{1}{2\sigma^2\,(2G(0)-1)} \, .$$
\end{theorem}
\smallskip

\begin{remark}
In this regime the deviation is entirely due to the moderate deviation behaviour of the scenery,
whereas the random walk does not contribute and behaves in a typical way. Asselah and
Castell~\cite{AC06} show that the regime in this result is maximal possible under Cram\'er's
condition, more precisely, higher regularity features of the scenery distribution
decide whether this behaviour persists when $b_n$ grows faster than $n^{2/3}$.
\end{remark}
\smallskip

\begin{remark}
For the sequence $b_n=n^\beta$ with $1/2 < \beta \le 2/3$, the
deviation speed $n^{2\beta-1}$, but not the rate function, in this result
was identified by Asselah and Castell~\cite{AC06} in $d\ge 5$ and by Asselah~\cite{As06} in $d=3$,
under the additional assumptions that the law of $\xi(0)$ has a symmetric density which
is decreasing on the positive half-axis.
\end{remark}

Turning to $d=2$, we define $\varkappa$ to be the optimal constant in the \emph{Gagliardo-Nirenberg} inequality,
$$\varkappa := \inf\big\{ c \, : \, \|f\|_4 \le c \, \|\nabla f\|_2^{\frac 12}\,\,  \|f\|_2^{\frac 12}
\, \mbox{ for all } f\in C^1_{\rm c}(\R^2) \big\}.$$
This constant features prominently in large deviation results for intersection local times
of Brownian motion and random walk intersection local times, see \cite{Ch04} for further discussion
of the Gagliardo-Nirenberg inequality and the associated constant~$\varkappa$.
\pagebreak[2]

\begin{theorem}[Moderate deviations in dimension $d=2$]\label{maind=2}\ \\ \vspace{-.3cm}
\begin{enumerate}
\item[(a)] If $n^{\frac 12}\sqrt{\log n} \ll b_n\ll n^{\frac 12}\log n$, then, as $n\uparrow\infty$,
$$\log\mathbb{P}\{X_n\geq b_n\}\sim-\frac{b_n^2}{n\log n} \, \frac{\pi (\det \Gamma)^{1/2}}{2\sigma^2}.$$
\item[(b)] If $n^{\frac 12}\log n\ll b_n\ll n/\log n$, then, as $n\uparrow\infty$,
$$\log\mathbb{P}\{X_n\geq b_n\}\sim-\frac{b_n}{\sqrt{n}}\, \frac{ (\det \Gamma)^{1/4}}{\varkappa^2\sigma}.$$
\item[(c)] Finally, for every $a>0$,
$$\log\mathbb{P}\{X_n\geq an^{\frac 12}\log n\} \sim -I(a)\log n,$$
where
$$I(a):=\left\{ \begin{array}{ll} \displaystyle
\frac{\pi a^2 (\det \Gamma)^{1/2}}{2\sigma^2}, & \displaystyle\mbox{ for } a\leq
\frac{\sigma}{\pi\varkappa^2(\det \Gamma)^{1/4}},\\[1mm]
\displaystyle\frac{a\, (\det \Gamma)^{1/4}}{\sigma\varkappa^2}-\frac{1}{2\pi\varkappa^4}, &\mbox{ for }
\displaystyle a\geq \frac{\sigma}{\pi\varkappa^2(\det \Gamma)^{1/4}}.
\end{array}\right.$$
\end{enumerate}
\end{theorem}
\smallskip

\begin{remark}
In regime~(a) the deviation is due to the moderate deviation
behaviour of the scenery only, but in regimes~(b) and~(c) there is
an additional contraction of the walks to achieve the moderate
deviation. There is only a very small gap between our moderate
deviation regime and the large deviation regime  studied in
\cite{GKS05}: Assuming that \emph{all} exponential moments of
$\xi(0)$ are finite and $b_n=an$, for some $a>0$, they obtain a
large deviation principle with speed~$n^{1/2}$ and a rate function
which is strongly dependent on the moment generating function of the
scenery variable.
\end{remark}

\begin{remark}
In the special case of simple random walk in \emph{Gaussian} scenery,
Theorem~\ref{maind=2}(a) is known from \cite{GKS05}.
\end{remark}

The regime $n^{\frac 12}\sqrt{\log n} \ll b_n\ll n/\log n$, which we consider
in Theorem~\ref{maind=2}, is \emph{maximal} for a moderate deviation principle using only Cram\'er's condition.
The following large deviation principle shows that for $b_n\gg n/\log n$
finer features of the scenery distribution (in this particular case the constant~$D$) enter into
the large deviation rate.

\pagebreak[4]
\begin{prop}[Special large deviations for $d=2$]\label{special}
Assume that, for some $D>0$,
\begin{equation}\label{P1}
\log\prob\big\{ \xi(0) > x \big\}  \sim - D\, x\, \qquad\mbox{ as $x\uparrow\infty$,}
\end{equation}
and suppose that $(b_n \log n)/n \to \infty$ and $\log b_n/\log n \to \beta\in [1,2)$.
Then, as $n\uparrow\infty$,\\[1mm]
\begin{equation}\label{P2}
\log \prob\big\{ X_n \ge b_n \big\}
\sim - \Big(\frac{b_n}{\log n}\Big)^{1/2}\, \Big(\frac{8\,K_2 D}{2-\beta} \Big)^{1/2} \, ,
\end{equation}
provided the underlying random walk is such that the
limit $K_2:=\lim_{n\to\infty} \frac{\me[\ell_n(0)]}{\log n}\in(0,\infty)$
exists.
\end{prop}

\begin{remark}
Note that this result is the planar case of the regime
$$\log \prob\{ \xi(0) > x\} \sim  -D x^{\frac d2} \qquad\mbox{ as }x\uparrow\infty,$$
which is described as `delicate' in \cite[Remark 1.2]{GKS05}.
The proof of Proposition~\ref{special} is 
based on large deviation results for the maximum of the local times obtained in~\cite{GHK06}.
\end{remark}
\smallskip

The remainder of the paper is structured as follows. Section~\ref{s3} is devoted to
statements about self-intersection local times of our random walk, which are of independent
interest. The proofs of our three theorems and Proposition~\ref{special}
follow in the subsequent four sections.
\smallskip

Throughout this paper we use the symbols $\mathsf P$ and $\mathsf E$
to denote probabilities, resp.~expectations, with respect to the
scenery variables only, and the symbols $\mathbb P$ and $\mathbb E$
to denote probabilities, resp.~expectations, with respect to both
the random walk and scenery.

We use the letters $c, C$ to denote positive, finite constants,
whose value can change at every occurrence, and which never depend
on random quantities. For nonnegative functions $f_n$, $g_n$,
possibly depending on the sampled walk or scenery, the Landau
symbols $f_n=o(g_n)$ and $f_n=O(g_n)$ denote $\lim f_n/g_n=0$,
respectively $\limsup f_n/g_n<\infty$, \emph{uniformly in the
sampled walk or scenery}. \pagebreak[1]

\section{Concentration inequalities for self-intersection local times}\label{s3}

Recall that $\{S_n \colon n\ge 0\}$ is a symmetric, aperiodic
random walk on the lattice~$\Z^d$,~$d\ge 2$, with nondegenerate
covariance matrix~$\Gamma$.  For integers $q>1$ we define the
\emph{$q$-fold self-intersection local time}
$\{\ell_n^{\ssup q} \colon n\ge 0\}$ of the random walk as
$$\ell_n^{\ssup q} := \sum_{z\in\Z^d} \ell^q_n(z)
= \sum_{1\le i_1, \cdots,i_q \le n} \1\big\{S_{i_1}=\cdots = S_{i_q}\big\}
\qquad \mbox{ for } n\ge 0\, .$$
We also denote the \emph{maximum} of the local times by
$$\ell_n^{\ssup \infty} := \max_{z\in\Z^d} \ell_n(z)\, .$$
The most important quantity is $\{\ell_n^{\ssup 2} \colon n\ge 0\}$,
which is simply called the self-intersection local time. Its asymptotic
expectations are
\begin{equation}\label{expd}
\mathbb{E} \ell_n^{\ssup 2} \sim
\left\{\begin{array}{ll}  n \, (2G(0)-1) & \mbox{ if } d\ge 3\, ,\\[2mm]
n\log n  \, \frac{1}{\pi\sqrt{\det\Gamma}} & \mbox{ if } d=2\, .
\end{array}\right.
\end{equation}
In $d\ge 3$ this is easy, for $d=2$ in the strongly aperiodic case
this follows from the local central limit theorem in the form
$\prob\{ S_n=0\} = 1/(n\, 2\pi\sqrt{\det \Gamma})  + o(1/n),$
see~\cite[Proposition~P7.9, p.75]{Sp76}, and can be extended to the
periodic case using Spitzer's trick, see~\cite[proof of
Proposition~P26.1, p.310]{Sp76}.
\smallskip

The main results of this section are the following concentration
inequalities for double and triple self-intersection local times,
which are of independent interest. They are therefore given in
somewhat greater generality than needed for
the proof of our main results. 
\pagebreak[4]

\begin{prop}[Concentration inequalities]\label{centredmoment}
Let $n\ge 2$. There exists a constant $c>0$ such that,
\begin{itemize}
\item[(a)] if $d>4$, then for $x\ge n^{\frac 23}\log^2n$,
$$\prob\big\{  |\ell_n^{\ssup 2} - \me\ell_n^{\ssup 2} |  \ge x \big\}
\le \exp\Big\{ -c\, \frac{x^{\frac 12}}{\log n} \Big\}\, ;$$
\item[(b)] if $d=4$, then for $x\ge n^{\frac 23}\log^3n$,
$$\prob\big\{  |\ell_n^{\ssup 2} - \me\ell_n^{\ssup 2}|  \ge x \big\}
\le \exp\Big\{ -c\, \frac{x^{\frac 12}}{\log^{3/2} n} \Big\}\, ;$$
\item[(c)] if $d=3$, then for $x\ge n^{\frac 12} \log^{9/2} n$,
$$\prob\big\{  |\ell_n^{\ssup 2} - \me\ell_n^{\ssup 2}|  \ge x \big\}
\le \exp\Big\{ -c\, \frac{x^{\frac 23}}{n^{\ui{\frac 13}}} \Big\}\, ; $$
\item[(d)] if $d> 4$, then for $x\ge n^{\frac 35}\log^2 n$,
$$\prob\big\{  |\ell_n^{\ssup 3} - \me\ell_n^{\ssup 3}|  \ge x \big\}
\le \exp\Big\{ -c\, \frac{x^{\frac 13}}{\log^{2/3} n} \Big\}\, ; $$ 
\item[(e)] if $d= 4$, then for $x\ge n^{\frac 35}\log^{7/2} n$,
$$\prob\big\{  |\ell_n^{\ssup 3} - \me\ell_n^{\ssup 3}|  \ge x \big\}
\le \exp\Big\{ -c\, \frac{x^{\frac 13}}{\log^{7/6} n} \Big\}\, . $$
\end{itemize}
\end{prop}


\begin{remark}
All of these inequalities are, to the best of our knowledge, new.
Similar concentration inequalities, but only for simple random walk and under considerably
stronger assumptions on the relationship of $x$ and $n$, have been found by
Asselah and Castell in~\cite[Propositions~1.4 and 1.6]{AC06} if $d\ge 5$,
and by Asselah in~\cite[Proposition 1.1]{As06} if $d=3$. In particular,
if $d\ge 5$, for the special case $x=yn$ they obtain an upper bound of
$\exp\{-c\sqrt{n}\}$, which is an improvement of (a). The proofs in
\cite{As06, AC06} are based on a delicate and powerful analysis of the number of sites in $\Z^d$
visited a certain number of times, and are therefore of independent interest.
In this paper we give a direct proof of Proposition~\ref{centredmoment},
which entirely avoids the discussion of the number of visits to individual sites,
and is therefore much easier than the method of Asselah and Castell.
\end{remark}

\subsection{Proof of Proposition~\ref{centredmoment}}

We start with some useful estimates for the partial Green's functions,
$$G_n(x):=\sum_{k=0}^n \prob\{ S_k=x\}, \qquad \mbox{ for } n\ge 2 \mbox{ and } x\in\Z^d .$$
\begin{lemma}\label{help}
For all $n\ge 2$,
$$\sum_{z\in\Z^d} G^2_n(z) \le \left\{ \begin{array}{ll}
C \, \sqrt{n} & \mbox{ if } d=3,\, \\
C\, \log n & \mbox{ if } d=4,\, \\
C\, & \mbox{ if } d>4.\end{array}\right. $$
\end{lemma}

\begin{Proof}{Proof.}
If $d=3$ we have from \cite[Proposition~P26.1, p.308]{Sp76} that
$G(z) \le C/(1+|z|)$. Then
$$\sum_{z\in\Z^3} G_n^2(z) = \sum_{|z|\le\sqrt{n}} G_n^2(z) + \sum_{|z|>\sqrt{n}} G_n^2(z)
 \le \sum_{|z|\le\sqrt{n}} G^2(z) + \Big(\sup_{|z|>\sqrt{n}} G(z) \Big) \,  \sum_{|z|>\sqrt{n}} G_n(z)\, .$$
The estimate for $G(z)$ shows that the first sum on the right is bounded by $C \sqrt{n}$. We
further have, from the definition of~$G_n$ and Chebyshev's inequality,
$$\begin{aligned}
\Big( \sup_{|z|>\sqrt{n}} G(z) \Big)\,  \sum_{|z|>\sqrt{n}} G_n(z) & \le C \, n^{-1/2}\,
 \sum_{k=0}^n \prob\{ |S_k| > \sqrt{n} \} \le C \,  n^{-1/2}\, \sum_{k=0}^n \frac{\me|S_k|^2}{n}
 \le C \, \sqrt{n}\, ,\end{aligned}$$
which completes the argument. In dimension $d\ge 4$ we use that, by~\cite[(1.4)]{Uc98}, we have
\begin{equation}\label{G}
G(z) \le \sum_{x\in\Z^d} \frac{\pi(x)}{1+|x-z|^{d-2}}\qquad \mbox{
for all }z\in\Z^d, \end{equation}%
where $(\pi(x) \colon x\in\Z^d)$ is a summable family of nonnegative
weights. If $d>4$, by the triangle inequality,
$$\Big(\sum_{z\in\Z^d} G^2(z) \Big)^{1/2} \le \sum_{x\in\Z^d}
\Big( \sum_{z\in\Z^d} \frac{\pi^2(x)}{(1+|x-z|^{d-2})^2} \Big)^{1/2}
= \Big( \sum_{x\in\Z^d} \pi(x)\Big)\,   \Big( \sum_{z\in\Z^d} \frac{1}{(1+|z|^{d-2})^2} \Big)^{1/2},$$
which is bounded by a constant. If $d=4$ we use first that
$$\sum_{z\in\Z^4} G_n^2(z) = \sum_{|z|\le n} G_n^2(z) + \sum_{|z|>n} G_n^2(z)
 \le \sum_{|z|\le n} G^2(z) + \Big(\sup_{z\in\Z^4} G(z) \Big) \,  \sum_{|z|>n} G_n(z)\, .$$
Clearly, $G$ is bounded, see (\ref{G}), and an argument analogous to
the case $d=3$ shows that the second sum on the right is bounded by
a constant. Using the triangle inequality as in the case $d>4$ we
obtain for the first sum on the right
$$\Big(\sum_{|z|\le n} G^2(z) \Big)^{1/2} \le \sum_{x\in\Z^4} \pi(x) \,
\Big(\sum_{|z+x|\le n} \frac{1}{(1+|z|^2)^2} \Big)^{1/2}\, .$$
It suffices to show that the content of the round bracket on the right is bounded
by a constant multiple of $\log n$, uniformly in~$x\in\Z^4$. On the
one hand, if $|x|\le 2n$ this follows easily from the fact that the sum can now
be taken over all $z\in\Z^4$ with $|z|\le 3n$. On the other hand, if $|x|> 2n$
the sum can be taken over the annulus $|x|-n \le |z| \le |x|+n$ and is thus easily
seen to be bounded by a constant.
\qed\end{Proof}

The proof of Proposition~\ref{centredmoment} requires the following `folklore' lemma
about the intersection of two independent random walks $\{S_n \colon n\ge 0\}$
and $\{S'_n \colon n\ge 0\}$ with $S_0=S_0'$. Denote
$$A_n := \sum_{i=1}^n\sum_{j=0}^{n-1} \1\{ S_i=S_j'\}\quad \mbox{ for } n\ge 1.$$
\pagebreak[2]
\begin{lemma}\label{moment}
There exists a constant $\vartheta>0$ such that,
\begin{itemize}
\item[(a)] if $d>4$, then
$\displaystyle\sup_{n\ge 2} \me \exp\big\{ \vartheta A_n^{1/2} \big\} <\infty\, ;$\\[1mm]
\item[(b)] if $d=4$, then
$\displaystyle \sup_{n\ge 2} \me \exp\big\{ \vartheta \sfrac{1}{\sqrt{\log n}} \,A_n^{1/2} \big\} <\infty\, ;$\\[1mm]
\item[(c)] if $d=3$, then
$\displaystyle \sup_{n\ge 2} \me \exp\big\{ \vartheta \, \big(\sfrac{A_n}{\sqrt n}\big)^{2/3} \big\} <\infty.$
\end{itemize}
\end{lemma}

\begin{Proof}{Proof.}
{F}rom the definition of $A_n$ we obtain, for moments of order $m\ge 1$,
$$\begin{aligned}
\me A_n^m &\le m! \sum_{1\le j_1 \le \cdots \le j_m\le n}\,
\sum_{0\le k_1,\ldots, k_m<n} \me \prod_{l=1}^m \1\{S_{j_l}=S'_{k_l}\} \\
&\le m! \sum_{\sigma\in\mathfrak{S}_m} \sum_{1\le j_1 \le \cdots \le j_m\le n}\,
\sum_{0\le k_1\le\cdots\le k_m< n} \, \sum_{x_1,\ldots, x_m} \me \prod_{l=1}^m \1\{S_{j_l}=x_l\}
\, \me \prod_{l=1}^m \1\{S'_{k_l}=x_{\sigma(l)}\} \\
& \le m! \sum_{\sigma\in\mathfrak{S}_m}
\sum_{x_1,\ldots, x_m} \prod_{l=1}^m G_n(x_l-x_{l-1}) G_n(x_{\sigma(l)}-x_{\sigma(l-1)}) ,\\
\end{aligned}$$
where $\mathfrak{S}_m$ denotes the group of all permutations of
$\{1,\ldots,m\}$, and we set $x_0:=0=:x_{\sigma(0)}$ for
convenience. Applying H\"older's inequality,
$$\me A_n^m \le (m!)^2 \, \sum_{x_1,\ldots, x_m} \prod_{l=1}^m G_n^2(x_l-x_{l-1}) =
(m!)^2 \, \Big( \sum_{x\in\Z^d} G_n^2(x)\Big)^m ,$$
and from Lemma~\ref{help} we obtain, for all~$n\ge 2$,
$$\me A_n^m \le \left\{ \begin{array}{ll}
(m!)^2\, C^m\, n^{m/2} & \mbox{ if } d=3,\, \\
(m!)^2\,C^m\,(\log n)^m & \mbox{ if } d=4,\, \\
(m!)^2\,C^m\, & \mbox{ if } d>4.\end{array}\right. $$
If $d>4$ this implies $\me\big(\sqrt{A_n}\big)^m \le \sqrt{ \me A_n^m } \le m!\, C^m,$
and~(a) follows by considering the exponential series. The analogous argument for $d=4$ gives~(b).
In $d=3$ we need an extra argument to complete the proof: We write
$\ell(m,n):=\lceil n/m \rceil +1$. Using an inequality of Chen,
\cite[Theorem~5.1]{Ch04} (with $p=2$ and $a=m$), we get, for $n\ge m$,
$$\begin{aligned}
\sqrt{\me A_n^m} & \le \sum_{\heap{k_1+\cdots +k_m=m}{k_1,\ldots,k_m\ge 0}}
\frac{m!}{k_1! \cdots k_m!} \sqrt{\me A^{k_1}_{\ell(m,n)} } \cdots \sqrt{\me A^{k_m}_{\ell(m,n)} } \\
& \le \sum_{\heap{k_1+\cdots +k_m=m}{k_1,\ldots,k_m\ge 0}}
\frac{m!}{k_1! \cdots k_m!} \sqrt{ (k_1!)^2\, C^{k_1} \ell(m,n)^{k_1/2}}\cdots
\sqrt{ (k_m!)^2\, C^{k_m} \ell(m,n)^{k_m/2}} \\
& \le \left(\mbox{$\heap{2m-1}{m}$}\right)\, m!\, C^m \, \big(\sfrac nm \big)^{m/4}
\le (m!)^{3/4}\, C^m\, n^{m/4},
\end{aligned}$$
and therefore $\me A_n^m \le (m!)^{3/2}\, C^m\,n^{m/2}.$
For $n\le m$ we get the same estimate immediately {f}rom the trivial
inequality $A_n^m\le n^{2m}\le (m!)^{3/2}\, C^m\,n^{m/2}$. We thus obtain,
for all $n,m$, that $$\me\big( n^{-1/3} \, A_n^{2/3}\big)^m
= n^{-m/3}\, \me\big( A_n^m\big)^{2/3}   \le m!\, C^m,$$
and~(c) follows by taking the exponential series.
\qed\end{Proof}
\medskip

Introduce, for $n\ge 1$,
$$\Lambda_n := \sum_{i=1}^n  \sum_{j,k=0}^{n-1}\1\{ S_i=S_j'=S_k'\}
\qquad\mbox{ and } \qquad \Lambda_n^{*} := \sum_{i=0}^{n-1}
\sum_{j,k=1}^{n}\1\{ S_i=S_j'=S_k'\}.$$ \pagebreak[3]

\begin{lemma}\label{moment3}
There exists a constant $\vartheta>0$ such that,
\begin{itemize}
\item[(a)] if $d>4$, then
$\displaystyle\sup_{n\ge 2} \me \exp\big\{ \vartheta \, \Lambda_n^{1/3} \big\} <\infty\,$;\\[1mm]
\item[(b)] if $d=4$, then
$\displaystyle \sup_{n \ge 2} \me \exp\big\{ \vartheta \, \sfrac{\Lambda_n^{1/3}}{(\log n)^{1/2}} \big\} <\infty$.
\end{itemize}
The same statements hold when $\Lambda_n$ is replaced by $\Lambda_n^*$.
\end{lemma}

\begin{Proof}{Proof.} We only consider $\Lambda_n$, as $\Lambda_n^*$ can be treated
analogously. {F}rom the definition of $\Lambda_n$ we obtain, for
moments of order~$m\ge 1$,
$$\begin{aligned}
\me \Lambda_n^m &\le m! \sum_{1\le j_1 \le \cdots \le j_m\le n}
\sum_{\heap{0\le k_1,\ldots, k_m<n}{0\le l_1,\ldots, l_m< n}}
\me \prod_{i=1}^m \1\{S_{j_i}=S'_{k_i}= S'_{l_i}\} \\
& \le m! \sum_{x_1,\ldots, x_m} \sum_{\heap{0\le k_1,\ldots, k_m< n}{0\le l_1,\ldots, l_m< n}}
\prod_{i=1}^m G_n(x_i-x_{i-1}) \,
\me \prod_{i=1}^m \1\{ S'_{k_i}= S'_{l_i} = x_i \}  ,\\
\end{aligned}$$
where we set $x_0:=0$ for convenience. Continuing with Cauchy-Schwarz, we get
$$\begin{aligned}
\le  m! \, \bigg( \sum_{x_1,\ldots, x_m} \prod_{i=1}^m G^2_n(x_i-x_{i-1}) \bigg)^{1/2} \,
\bigg( \sum_{x_1,\ldots, x_m} \Big( \sum_{\heap{0\le k_1,\ldots, k_m< n}{0\le l_1,\ldots, l_m<n}}
\me \prod_{i=1}^m \1\{ S'_{k_i}= S'_{l_i} = x_i \} \Big)^{2} \bigg)^{1/2}\, .
\end{aligned}$$
By~Lemma~\ref{help} the first bracket is bounded by $C^m$ if $d>4$, and
by $C^m (\log n)^m$ if $d=4$. To analyse the second bracket we denote
by $\mathcal T_m$ the set of all mappings $\tau\colon \{1,\ldots, 2m\}
\to \{1,\ldots, m\}$ such that $\#\tau^{-1}\{j\}=2$ for all $j\in\{1,\ldots,m\}$. 
For the cardinality of $\mathcal T_m$ we get
\begin{equation}\label{Tcount}
\# \mathcal T_m \le \left(\heap{2m}{m} \right) (m!)^2 \le C^m\,
(m!)^2\, .
\end{equation}
Given $(k_1,\ldots,k_m)$ and $(l_1,\ldots,l_m)$ there exists at least one
ordered tuple $(k'_1,\ldots, k'_{2m})$ with $k_1'\le \cdots \le k_{2m}'$ 
with $\{k_1,\ldots, k_m, l_1,\ldots, l_m\}=\{k'_1,\ldots, k'_{2m}\}$ and
$\tau\in\mathcal T_m$ such that $\tau(i)=j$ if $k'_i=l_j$ or $k'_i=k_j$. 
Hence we obtain,
$$\begin{aligned}
\sum_{\heap{0\le k_1,\ldots, k_m< n}{0\le l_1,\ldots, l_m< n}}
\me \prod_{i=1}^m \1\{ S'_{k_i}= S'_{l_i} = x_i \}
& \le\sum_{\tau\in\mathcal T_m} \sum_{0\le k'_1\le\cdots\le k'_{2m}< n} \,\,
\prod_{i=1}^{2m} \prob\big\{ S'_{k'_i}-S'_{k'_{i-1}}=x_{\tau(i)}-x_{\tau(i-1)} \big\}\\
& \le  \sum_{\tau\in\mathcal T_m}
\prod_{i=1}^{2m} G_n(x_{\tau(i)}-x_{\tau(i-1)}),
\end{aligned} $$
and, using the triangle inequality,
$$\begin{aligned}
\bigg( \sum_{x_1,\ldots, x_m}  \Big( \sum_{\tau\in\mathcal T_m} &
\prod_{i=1}^{2m} G_n(x_{\tau(i)}-x_{\tau(i-1)}) \Big)^2  \bigg)^{1/2}
\le \sum_{\tau\in\mathcal T_m} \bigg( \sum_{x_1,\ldots, x_m}
\prod_{i=1}^{2m} G^2_n(x_{\tau(i)}-x_{\tau(i-1)}) \bigg)^{1/2}\\
& \le  \# \mathcal T_m \, \bigg( \sum_{x_1,\ldots, x_{2m} }
\prod_{i=1}^{2m} G^2_n(x_{i}-x_{i-1}) \bigg)^{1/2}\, .
\end{aligned}$$
By Lemma~\ref{help} the bracket is 
bounded by $C^m$ if $d>4$, and 
by $C^m (\log n)^{2m}$ if $d=4$. Thus, together with
\eqref{Tcount}, we obtain the estimates
$$\me \Lambda_n^m  \le \left\{
\begin{array}{ll}
\big(m! \, C^m \big)^3 & \mbox{ if }\, d>4, \\[2pt]
\big(m! \, C^m \, (\log n)^{m/2}\big)^3 & \mbox{ if }\, d=4.
\end{array}  \right. $$
But $\me (\Lambda_n^{1/3})^m  \le \big(\me \Lambda_n^m\big)^{1/3}$,
and both statements follow by taking exponential series.
\qed\end{Proof}

For any $N\ge 0$ we use the classical decomposition
$$\ell_{2^N}^{\ssup 2} - \me  \ell_{2^N}^{\ssup 2}
= 2\, \sum_{j=1}^N \sum_{k=1}^{2^{j-1}} \overline{A}_{j,k},$$
where
$$\overline{A}_{j,k} := \overline{A}_{j,k}(N) :=
\sum_{\heap{(2k-2)2^{N-j}<l\le (2k-1)2^{N-j}}{ (2k-1)2^{N-j}<m \le (2k)2^{N-j}}}
\Big( \1\{S_l=S_m\} - \prob\{S_l=S_m\} \Big) .$$
For fixed $1\le j\le N$ the random variables
$\overline{A}_{j,k}$, for $k=1,\ldots, 2^{j-1}$, are independent, identically distributed
with the law of $A_{2^{N-j}}- \me A_{2^{N-j}}$. The next proposition exploits this independence,
and the moment results of Lemma~\ref{moment} to give large deviation upper bounds.
\pagebreak[2]

\begin{prop}[Large deviation upper bounds]\label{Nag}
For every $\eps>0$ there exists $c=c(\eps)>0$ such that,
for all $1\le j \le N$,
\begin{itemize}
\item[(a)] if $d>4$, then
$\,\,\displaystyle\prob\Big\{ \Big|\sum_{k=1}^{2^{j-1}} \overline{A}_{j,k}(N)\Big| \ge \eps \, x \Big\}
\le \exp\big\{ -c\, \sqrt{x} \big\}$ for all $x\ge (2^N)^{2/3}$;\\[1mm]
\item[(b)] if $d=4$, then
$\,\,\displaystyle \prob\Big\{ \Big|\sum_{k=1}^{2^{j-1}} \overline{A}_{j,k}(N)\Big| \ge \eps \, x\Big\}
\le \exp\Big\{ -c\, \sqrt{\frac{x}{N}} \Big\}$ for all $x\ge N(2^N)^{2/3}$;\\[1mm]
\item[(c)] if $d=3$, then
$\,\,\displaystyle \prob\Big\{ \Big|\sum_{k=1}^{2^{j-1}} \overline{A}_{j,k}(N) \Big|\ge \eps \, x\Big\}
\le \exp\Big\{ -c\, \frac{x^{2}}{2^{N}} \Big\} + \exp\Big\{ -c\, \frac{x^{2/3}2^{j/3}}{2^{N/3}} \Big\}$
\newline\hspace*{10cm} for all $x\ge N^{9/2}\,(2^{N})^{1/2}$.
\end{itemize}
\end{prop}

The proof of this result will be postponed to the next section.

\begin{Proof}{Completion of the proof of Proposition~\ref{centredmoment}(a)\,--\,(c).}
We use two simple ingredients, stated below as~\eqref{obs} and~\eqref{weis}.
\emph{First}, note that, for any $N\ge 0$ and any
choice of nonnegative  weights~$p_j$, $1\le j \le N$, with~$\sum p_j\le 1$,
we have
\begin{equation}\label{obs}\begin{aligned}
\prob\big\{ |\ell_{2^N}^{\ssup 2}-\me \ell_{2^N}^{\ssup 2}| \ge \eps y \big\} & =
\prob\Big\{ 2\Big| \sum_{j=1}^N \sum_{k=1}^{2^{j-1}} \overline{A}_{j,k} \Big| \ge \eps y \Big\}
\le \sum_{j=1}^N \prob\Big\{ \Big| \sum_{k=1}^{2^{j-1}} \overline{A}_{j,k} \Big| \ge \frac{\eps y p_j}{2} \Big\}.\\
\end{aligned}\end{equation}
\emph{Second}, for any $n\ge 2$ there exists the representation
$$n=2^{N_1} + \cdots + 2^{N_l},$$
where $l\ge 1$ and $N_1> \cdots >N_l \ge 0$ are integers. Note that $l\le c \log n$.
Write $n_0:=0$ and $n_i:=2^{N_1} + \cdots + 2^{N_i}$ for $1\le i \le l$, and denote
$$B_i:= \sum_{n_{i-1}<j<k\le n_i} \1\{S_j=S_k\}, \quad \mbox{ and } \quad
D_i:= \sum_{\heap{n_{i-1}<j\le n_i}{n_i < k \le n}} \1\{S_j=S_k\}.$$
Then $\sum_{1\le j<k\le n} \1\{S_j=S_k\} = \sum_{i=1}^l B_i + \sum_{i=1}^{l-1} D_i.$
We thus have, for any choice of nonnegative weights~$q_i$, $1\le i \le l$, with $\sum q_i\le 1$,
for $x$ large enough to satisfy $xq_i> 4\me D_i$,
\begin{equation}\label{weis}
\begin{aligned}
\prob\big\{ & |\ell_{n}^{\ssup 2}-\me \ell_{n}^{\ssup 2}| \ge  x \big\} \le
\sum_{i=1}^l \prob\big\{  |B_i-\me B_i| \ge \sfrac{x q_i}{4}  \big\}
+ \sum_{i=1}^{l-1} \prob\big\{  D_i\ge \sfrac{x q_i}{4} \big\} .
\end{aligned}
\end{equation}
Depending on the dimension, we use the ingredients \eqref{obs} and \eqref{weis} with
different choice of weights.  If $d=3$ we define $q_i=b 2^{(N_i-N_1)/2}$ with
$b=(\sum_{j=1}^\infty 2^{-j/2})^{-1}$, and apply~\eqref{obs} for
$$N=N_i, \,\,\, y=\frac{xq_i}{4\eps}  \mbox{ and weights }
p_j=aj^{-2}  \mbox{ with }  a=\Big(\sum_{j=1}^\infty j^{-2}\Big)^{-1},$$
where $\eps>0$ may be chosen independently of~$i, j$ such that
$yp_j/2\ge N_i^{9/2} (2^{N_i})^{1/2}$. Using \eqref{obs}, Proposition~\ref{Nag}~(c)
and that $l\le c \log n$, this gives
\begin{equation}\label{1sum}
\begin{aligned}
\sum_{i=1}^l \prob\big\{  |B_i-\me B_i| \ge  \sfrac{x q_i}{4}\big\} &
\le  \sum_{i=1}^l \sum_{j=1}^{N_i}  \exp\Big\{ -c\, \sfrac{(yp_j)^{2}}{2^{N_i}} \Big\}
+ \exp\Big\{ -c\, \sfrac{(yp_j)^{2/3}2^{j/3}}{2^{N_i/3}} \Big\} \\
& \le  \exp\Big\{ -c\, \sfrac{x^{2/3}}{n^{1/3}} \Big\} .
\end{aligned}
\end{equation}
As (with $\stackrel{d}{=}$ denoting equality of distributions)
$$D_i\stackrel{d}{=} \sum_{j=1}^{2^{N_i}} \sum_{k=1}^{n-n_i} \1\{ S_j=S'_k\}
\le \sum_{j=1}^{2^{N_i}} \sum_{k=0}^{2^{N_i}-1} \1\{ S_j=S'_k\} = A_{2^{N_i}},$$
the second sum in \eqref{weis} can be estimated using Chebyshev's inequality and Lemma~\ref{moment},
\begin{equation}\label{2sum}
\begin{aligned}
\sum_{i=1}^{l-1} \prob\big\{  D_i\ge \sfrac{x q_i}{4} \big\}
& \le \sum_{i=1}^{l-1} \prob\big\{ \sfrac{A_{2^{N_i}}}{2^{N_i/2}} \ge \sfrac{x q_i}{4 2^{N_i/2}} \big\} \\
& \le \sum_{i=1}^{l-1}\exp\big\{ -c \, \big( \sfrac{xq_i}{2^{N_i/2}}\big)^{2/3} \big\}
\le \exp \big\{ - c\, \sfrac{x^{2/3}}{n^{1/3}} \big\},
\end{aligned}\end{equation}
and the proof of (c) follows by plugging~\eqref{1sum} and~\eqref{2sum}
into~\eqref{weis}.  The proof of (a), (b) is analogous, but now the weights
are chosen to be equal, i.e. $p_j=1/N$ and  $q_i=1/l$. We leave the obvious details
to the reader.
\qed\end{Proof}
\medskip

An analogous argument can be carried out for triple self-intersections. Indeed,
for any $N\ge 0$ we have
\begin{equation}\label{zzz}
\ell_{2^N}^{\ssup 3} - \me \ell_{2^N}^{\ssup 3} = \sum_{j=1}^N \sum_{k=1}^{2^{j-1}} \overline{\Lambda}_{j,k}
+ \sum_{j=1}^N \sum_{k=1}^{2^{j-1}} \overline{\Lambda}^{*}_{j,k}
\end{equation}
where
$$\overline{\Lambda}_{j,k} := \sum_{\heap{(2k-2)2^{N-j}<l\le (2k-1)2^{N-j}}{ (2k-1)2^{N-j}<m,n \le (2k)2^{N-j}}}
\Big( \1\{S_l=S_m=S_n\} - \prob\{S_l=S_m=S_n\} \Big) $$
and
$$\overline{\Lambda}^{*}_{j,k} := \sum_{\heap{(2k-2)2^{N-j}<l,m\le (2k-1)2^{N-j}}{ (2k-1)2^{N-j}<n \le (2k)2^{N-j}}}
\Big( \1\{S_l=S_m=S_n\} - \prob\{S_l=S_m=S_n\} \Big) .$$
Again, for fixed $1\le j\le N$ the random variables
$\overline{\Lambda}_{j,k}$, for $k=1,\ldots, 2^{j-1}$, are independent, identically distributed
with the law of $\Lambda_{2^{N-j}}- \me \Lambda_{2^{N-j}}$, and
 the random variables
$\overline{\Lambda}^{*}_{j,k}$, for $k=1,\ldots, 2^{j-1}$, are independent, identically distributed
with the law of $\Lambda^*_{2^{N-j}}- \me \Lambda^*_{2^{N-j}}$.

\pagebreak[2]

\begin{prop}[Large deviation upper bounds]\label{Nag2}
For any $\eps>0$ there exists $c=c(\eps)>0$ such~that, for all
$1\le j\le N$,
\begin{itemize}
\item[(a)] if $d>4$, then
$\,\,\displaystyle\prob\Big\{ \Big|\sum_{k=1}^{2^{j-1}} \overline{\Lambda}_{j,k}\Big| \ge \eps \, x \Big\}
\le \exp\big\{ -c\,x^{1/3}  \big\} \, ,$ for all $x\ge (2^N)^{3/5}$;\\[1mm]
\item[(b)] if $d=4$, then
$\,\,\displaystyle \prob\Big\{ \Big|\sum_{k=1}^{2^{j-1}}
\overline{\Lambda}_{j,k}\Big| \ge \eps \, x\Big\} 
\le \exp\big\{ -c\, \big(\sfrac x{N^{3/2}}\big)^{1/3}  \big\} \, ,$ for all $x\ge N^{3/2} (2^N)^{3/5} $.\\
\end{itemize}
The same estimates hold for $\overline{\Lambda}_{j,k}$ replaced by $\overline{\Lambda}^{*}_{j,k}$.
\end{prop}

Again we postpone the proof of Proposition~\ref{Nag2} to the next section
and first complete the details of the remaining parts of Proposition~\ref{centredmoment}.

\begin{Proof}{Proof of Proposition~\ref{centredmoment}(d),(e).}
For any $N\ge 0$, we have by \eqref{zzz},
\begin{equation}\label{ubs}\begin{aligned}
\prob\big\{ |\ell_{2^N}^{\ssup 3}-\me \ell_{2^N}^{\ssup 3}| \ge \eps y \big\} & =
\prob\Big\{ \Big| \sum_{j=1}^N \sum_{k=1}^{2^{j-1}} \overline{\Lambda}_{j,k} \Big| \ge \frac{\eps y}2  \Big\}
+  \prob\Big\{ \Big| \sum_{j=1}^N \sum_{k=1}^{2^{j-1}} \overline{\Lambda}^*_{j,k} \Big| \ge \frac{\eps y}2  \Big\} \\
& \le \sum_{j=1}^N \prob\Big\{ \Big| \sum_{k=1}^{2^{j-1}} \overline{\Lambda}_{j,k} \Big| \ge \frac{\eps y}{2N} \Big\}
+ \sum_{j=1}^N \prob\Big\{ \Big| \sum_{k=1}^{2^{j-1}} \overline{\Lambda}^*_{j,k} \Big| \ge \frac{\eps y}{2N} \Big\}.
\end{aligned}\end{equation}
For any $n\ge 2$ there exists the representation
$n=2^{N_1} + \cdots + 2^{N_l},$ where $N_1> \cdots >N_l \ge 0$ are integers. Note that $l\le c \log n$.
Write $n_0:=0$ and $n_i:=2^{N_1} + \cdots + 2^{N_i}$ for $1\le i \le l$, and denote
$$B_i:= \sum_{n_{i-1}<j,k,l\le n_i} \1\{S_j=S_k=S_l\},$$
$$ D_i:= \sum_{\heap{n_{i-1}<j,k\le n_i}{n_i < l \le n}} \1\{S_j=S_k=S_l\} \qquad \mbox{ and } \qquad
E_i:= \sum_{\heap{n_{i-1}<j\le n_i}{n_i < k,l \le n}} \1\{S_j=S_k=S_l\}.$$
Then $\ell_n^{\ssup 3} = \sum_{i=1}^l B_i + \sum_{i=1}^{l-1} D_i
+ \sum_{i=1}^{l-1} E_i.$ As $\me D_i$ and $\me E_i$ are bounded by a
constant multiple of $\log n$,  we get for all sufficiently large $x$,
\begin{equation}\label{wo}
\begin{aligned}
\prob\big\{ & |\ell_{n}^{\ssup 3}-\me \ell_{n}^{\ssup 3}| \ge  x \big\} \le
\sum_{i=1}^l \prob\big\{  |B_i-\me B_i| \ge \sfrac{x}{3l}  \big\}
+ \sum_{i=1}^{l-1} \prob\big\{  D_i \ge \sfrac{x}{3l}  \big\}
+ \sum_{i=1}^{l-1} \prob\big\{  E_i\ge \sfrac{x}{3l} \big\} .
\end{aligned}
\end{equation}
We now look at the case $d=4$. Using \eqref{ubs} with $y=x/(3l\eps)$, Proposition~\ref{Nag2}(b)
and that $l\le c \log n$, this gives
\begin{equation}\label{01sum}
\begin{aligned}
\sum_{i=1}^l \prob\big\{  |B_i-\me B_i| \ge  \sfrac{x}{3l}\big\} &
\le  2\sum_{i=1}^l \sum_{j=1}^{N_i}  \exp\Big\{ -c\, 
\big(\frac{x}{lN_i^{5/2}}\big)^{1/3} \Big\} \le  \exp\Big\{ -c\,
\frac{x^{1/3}}{\log^{7/6} n} \Big\} .
\end{aligned}
\end{equation}
As we have
$$D_i\stackrel{d}{=}\, \sum_{j,k=0}^{2^{N_i}-1}\, \sum_{m=1}^{n-n_i} \1\{ S_j=S_k=S'_m\}
\le \sum_{j,k=0}^{2^{N_i}-1} \sum_{m=1}^{2^{N_i}} \1\{
S_j=S_k=S'_m\} = \Lambda^*_{2^{N_i}},$$
the second sum in \eqref{wo}
can be estimated using Chebyshev's inequality and
Lemma~\ref{moment3}(b),
\begin{equation}\label{02sum}
\begin{aligned}
\sum_{i=1}^{l-1} \prob\big\{  D_i\ge \sfrac{x}{3l} \big\} & \le 
\sum_{i=1}^{l-1} \prob\Big\{ \frac{\Lambda^*_{2^{N_i}}}{N_i^{3/2}}
\ge \frac{x}{3lN_i^{3/2}} \Big\} \le l\, \exp\Big\{ -c \, \Big(
\frac{x}{lN_1^{3/2}}\Big)^{1/3} \Big\} \le \exp \Big\{ - c\,
\frac{x^{1/3}}{\log ^{5/6} n} \Big\}.
\end{aligned}\end{equation}
The same estimate holds for $E_i$ in place of $D_i$, using the estimate for $\Lambda_{2^{N_i}}$
instead of  $\Lambda^*_{2^{N_i}}$. The proof of (c) follows by plugging this, \eqref{02sum} and~\eqref{01sum}
into~\eqref{wo}. The case $d\ge 5$ is analogous.\qed\end{Proof}

\subsection{Proof of Propositions~\ref{Nag} and~\ref{Nag2}}

\begin{Proof}{Proof of Proposition~\ref{Nag}.}
We first give the argument in the case $d\ge 5$. Take a continuously differentiable
function $g\colon(0,\infty)\to\R$ with non-increasing derivative, such that
\begin{itemize}
\item[(a)] $g'(x) > 2/x$ for all $x>0$,
\item[(b)] $g(x)= \vartheta \sqrt{x}$ for all $x\ge x_0$,
\end{itemize}
where $\vartheta$ is chosen as in Lemma~\ref{moment}. For
$1\le j \le N$ denote
$$b_j(N):=\me\Big[ \exp\big\{ g\big( \overline{A}_{j,1}(N)\big)\big\}
\, \1\{ \overline{A}_{j,1}(N) >0 \}\Big],$$
and recall from Lemma~\ref{moment}(a) that $b_j(N)$ is uniformly bounded in $j$ and $N$.
By Theorem 2.3 of \cite{Na79} (with $\gamma_1=\gamma_2=\gamma_3=1/3$,
$\gamma=2/3$ and $\delta=2$) we obtain the bound
\begin{align}
\prob\Big\{ \sum_{k=1}^{2^{j-1}} \overline{A}_{j,k} \ge \eps \, x\Big\}
& \le e^{1/2}\,\exp\Big\{ - \sfrac{a^2 \,\eps^2\, x^2}{2(a+1) 2^{j-1} V_j(N)} \Big\}\label{1a}\\
&\qquad + e^{1/2}\,\exp\Big\{ - \sfrac{2a \, \eps \, x}{3 S^{-1}\big(\frac{a \eps x}{3 e^a 2^{j-1} b_j(N)}\big) } \Big\} \label{1b}\\
&\qquad + 2^j \,b_j(N)\, e^{1/2}\, \exp\Big\{ -g\big( \sfrac 23 \,\eps\,x\big)  \Big\}
 + 2^{j-1}\, \prob\Big\{ \overline{A}_{j,1} \ge \sfrac 23 \, \eps\, x \Big\},\label{1d}
\end{align}
where $V_j(N)$ is the variance of $\overline{A}_{j,1}$, the constant $a$ is the unique solution of the equation
$(u+1)=e ^{u-1}$, and $S^{-1}$ is the inverse of the strictly decreasing function $u \mapsto S(u):=e^{-g(u)} g'(u) u^2$,
see \cite[p.765]{Na79}. By Chebyshev's inequality,
$$\prob\big\{ \overline{A}_{j,1} \ge x \big\} \le \big( \sup_N \sup_{j\le N} b_j(N) \big) \, e^{-g(x)},$$
and therefore the two terms in \eqref{1d} are bounded by a constant multiple of
$$2^N\, \exp\Big\{ -g\big( \sfrac 23\, \eps \,x\big) \Big\} \qquad \mbox{ for all } j\le N.$$
Recalling the definition of $g$ we arrive at an upper bound of
\begin{equation}\label{1cd}
C \, \exp\big\{ -c\, \,\sqrt{x} \big\} \qquad \mbox{ for all } N\ge 1.
\end{equation}
If $x\ge (2^N)^{2/3}$, then $x^2/2^{j-1}=x^{1/2}\, x^{3/2}/2^{j-1} > \sqrt{x}$
for all $j\le N$. Further, using this inequality and the boundedness of $V_j(N)$,
the term in~\eqref{1a} is also bounded by a constant multiple of $\exp\{-c\sqrt{x}\}$.

To show that also the term in \eqref{1b} is negligible, recall that the function $S$
is strictly decreasing.  Hence, the term in \eqref{1b} is bounded by
$$C \, \exp\Big\{ - c \frac{x}{S^{-1}\big(\frac{c}{2^{N/3}} \big)} \Big\}.$$
{F}rom the definition of the functions $g$ and $S$ it is easy to see that
$$S^{-1}\big(\sfrac{c}{2^{N/3}} \big) \le C N^2.$$
This implies that the term in \eqref{1b} is bounded by a constant multiple
of $\exp\{ -c \, x/N^2\}$, and is therefore also negligible compared to \eqref{1cd}.
This completes the bound for $\sum \overline{A}_{j,k}$.
The same reasoning can be applied with $-\overline{A}_{j,k}$ in place of
$\overline{A}_{j,k}$, using only the trivial fact that $-\overline{A}_{j,1}$
is bounded from above, uniformly in $j$. Hence we get the same bound for
$-\sum \overline{A}_{j,k}$. This completes the proof in dimensions~\mbox{$d\ge5$}.
The result in $d=4$ is a modification of this argument,
using the random variable $(N-j)^{-1} \overline{A}_{j,k}$ instead of $\overline{A}_{j,k}$,
and details are left to the reader.

Turning to dimension~$d=3$, we use that
$$\prob\Big\{ \sum_{k=1}^{2^{j-1}} \overline{A}_{j,k} \ge \eps x \Big\}
= \prob\Big\{ \sum_{k=1}^{2^{j-1}} \frac{ \overline{A}_{j,k}} {2^{(N-j)/2}} \ge
\eps \frac{x}{2^{(N-j)/2}} \Big\},$$
and choose a function~$g\colon (0,\infty)\to\R$ which satisfies the same conditions
as above, except that we now replace condition~(b) by~$g(x)=\vartheta\,x^{2/3}$ for all
$x\ge x_0$, and $\vartheta$ as in Lemma~\ref{moment}. We define
$$b_j(N) := \me\big[ \exp\big\{g\big(\overline{A}_{j,1}/2^{(N-j)/2}\big)\,
 \1\{ \overline{A}_{j,1} >0 \} \big],$$
and by Theorem~\cite[Theorem 2.3]{Na79} we obtain
\begin{align}
\prob\Big\{ \sum_{k=1}^{2^{j-1}} \frac{\overline{A}_{j,k}}{2^{(N-j)/2}} \ge \eps \, & \,
\frac{x}{2^{(N-j)/2}}\Big\} \le \exp\Big\{ - c \, \frac{x^2}{2^N} \Big\}
+\exp\Big\{ - c \,\frac{x}{2^{(N-j)/2}\,  S^{-1}\big(\frac{c x}{2^{(N+j)/2}}\big) } \Big\} \label{01b}\\
& + C\,2^j \,b_j(N)\, \exp\Big\{ -g\big(c\, \frac{x}{2^{(N-j)/2}}\big)  \Big\}
 + 2^{j-1}\, \prob\Big\{ \frac{\overline{A}_{j,1}}{2^{(N-j)/2}}
\ge c\,\frac{x}{2^{(N-j)/2}} \Big\}.\label{01d}
\end{align}
The two terms in \eqref{01d} are bounded by $2^N\,\exp\{-c\,x^{2/3}/2^{(N-j)/3}\}$. To bound
the last term in \eqref{01b} we use that, for $x\ge 2^{N/2}/N^2$,
$$S^{-1}\Big( \frac{cx}{2^{(N+j)/2}}\Big) \le S^{-1}\Big( \frac{cx}{2^N}\Big) \le
S^{-1}\Big(\frac{c}{N^2 2^{N/2}}\Big) \le C N^{3/2}\, ,$$
to get
$$\exp\Big\{ - c \frac{x}{2^{(N-j)/2}\,  S^{-1}\big(\frac{c x}{2^{(N+j)/2}}\big) } \Big\}
\le \exp\Big\{ -c \, \frac{x 2^{j/2}}{2^{N/2} N^{3/2}} \Big\}\, .$$
As $x\ge 2^{N/2} N^{9/2}$ this term is also bounded by $\exp\{-c\,x^{2/3}/2^{(N-j)/3}\}$,
completing the proof. \qed\end{Proof}

\begin{Proof}{Proof of Proposition~\ref{Nag2}.} We use the same arguments as in
Proposition~\ref{Nag}, but now for a function $g\colon (0,\infty)\to\R$ with
condition~(b) replaced by~$g(x)=\vartheta x^{1/3}$ for $x\ge x_0$. Then both
terms in \eqref{1d} give contributions bounded by $\exp\{ -c\, x^{1/3}\}$. If
$x\ge (2^N)^{3/5}$, then $x^2/2^{j-1} \ge x^{1/3}$, and hence we obtain the
same bound for~\eqref{1a}. Under the same condition~$x\ge (2^N)^{3/5}$, we have
$$S^{-1}\big(cx/2^{j-1}\big)\le S^{-1}\big(c/(2^N)^{2/5}\big) \le C\, N^3\, ,$$
hence the term in \eqref{1b} is of smaller order.
\qed\end{Proof}

\subsection{A large deviation bound for the maximum of the local times}

We complete this section with an easy lemma, which provides bounds for the large
deviation probabilities of the
maximum $\ell_n^{\ssup \infty}$ of the local times. Ideas for this
proof are taken from Gantert and Zeitouni~\cite{GZ98}.
\pagebreak[2]
\begin{lemma}[Large deviation bounds for the maximal local time]\label{maxdev}
There exists $c>0$ such that
\begin{itemize}
\item[(a)] if $d\ge 3$, then for each sequence $a_n\to\infty$
and all $n\ge 2$,
$$\prob\big\{  \ell_n^{\ssup \infty}  > a_n \big\} \le n\,\exp\big\{ -c\, a_n \big\}\, ;$$\\[-6mm]
\item[(b)] if $d=2$, then for each sequence $a_n/\log n\to\infty$ and all $n\ge 2$,
$$\prob\big\{  \ell_n^{\ssup \infty}  > a_n \big\}
\le n\, \exp\Big\{ -c\, \frac{a_n}{\log n } \Big\}\, .$$
\end{itemize}
\end{lemma}

\begin{Proof}{Proof.}
Without loss of generality we may assume that all $a_n$ are positive integers.
We first reduce the problem to a large deviation bound for $\ell_n(0)$. Defining the
stopping times $T_z:=\min\{ k\ge 1 \colon S_k=z \}$ we have, for all nonnegative integers~$x$,
$$\begin{aligned}
\prob\big\{  \ell_n^{\ssup \infty}  >x \big\}
& \le \sum_{z\in\Z^d} \prob\big\{ \ell_n(z) > x \big\}
 = \sum_{z\in\Z^d} \sum_{k=1}^n \prob\{T_z=k\}\, \prob\big\{ \ell_{n-k}(0) \ge x \big\}\\
& \le  \prob\big\{ \ell_{n}(0) \ge x \big\} \, \sum_{z\in\Z^d} \prob\{T_z\le n\}  \, .
\end{aligned}$$
Now $\sum_z \prob\{T_z\le n\} \le \sum_z \sum_{k=1}^n \prob\{ S_k=z\} =n$, so
that it suffices to bound the large deviation probabilities of $\ell_n(0)$.
By the strong Markov property applied at the successive hitting times of the origin,
we get
\begin{equation}\label{queye}
\prob\big\{ \ell_{n}(0) \ge a_n \big\} \le  \prob\big\{ T_0 \le n \big\}^{a_n} \, .
\end{equation}
In the transient case, $d\ge 3$, this gives~(a) with $c:=-\log \prob\{ T_0<\infty\}>0$.
In the recurrent case $d=2$, we use the last exit decomposition, for all $2\le k \le n$,
$$1\le \sum_{j=0}^k \prob\{ S_j=0\}\, \prob\{ \ell_{n-k}(0) =0 \}
+ \sum_{j=k+1}^{n} \prob\{ S_j=0\}\, .$$ By \cite[Proposition~P7.6,
p.72]{Sp76} we have $\prob\{ S_j=0\}\le \frac cj$ for $j\ge 1$. This
implies that
$$(\log k) \, \prob\{ \ell_{n-k}(0)=0 \} \ge C \Big[ 1 - c\, \Big( \sum_{j=k+1}^{n} \frac 1j \Big) \Big]\, .$$
Now let $k=\lceil \eta n\rceil$ and choose $\eta\in(0,1)$ sufficiently close to one, so that
the right hand side is bounded from zero by a positive constant. Hence,
$$\prob\{ T_0> n(1-\eta) \} =\prob\{ \ell_{\lfloor n(1-\eta)\rfloor}(0)=0 \}
\ge \frac{c}{\log n}\, ,$$
and thus $\log \prob\{ T_0 \le n\} = \log (1- \prob\{ T_0> n\})
\le -c/\log n$. Plugging this into \eqref{queye} completes the proof of~(b).
\qed\end{Proof}
\smallskip

\section{Precise asymptotics in dimensions $d\ge 4$: Proof of Theorem~\ref{maind>3}}

The main ingredient of the proof is the following proposition. Recall that the probability~$\mathsf P$
refers exclusively to the scenery variables with fixed random walk samples, and the Landau symbols
are uniform in these samples.

\begin{prop}\label{vitali}
Assume that, for some $A>0$ and all sufficiently large~$n$,
$$\Gamma_n:=\sum_{z\in\Z^d} \ell_n^3(z) \le n \log^2 n\qquad \mbox{ and } \qquad
V_n^2:=\sigma^2 \sum_{z\in\Z^d} \ell_n^2(z)\le An\, .$$ Then, for $\sqrt{n} \ll b_n \ll n^{2/3}/\log^{3/2} n$,
we have
\begin{equation}\label{cramer}
{\mathsf P}\Big\{ \sum_{z\in\Z^d} \ell_n(z) \xi(z) \ge b_n \Big\}
=\frac{V_n}{\sqrt{2\pi} b_n} \, \exp\Big\{ - \frac{b_n^2}{2V_n^2} \Big\} \, (1+o(1))\, .
\end{equation}
\end{prop}

\begin{Proof}{Proof of Theorem~\ref{maind>3}.}
On the event $$\big\{ |\ell_n^{\ssup 2}- \me \ell_n^{\ssup 2}| \le n^{2/3} \log^3 n, \,
\ell_n^{\ssup 3} \le n\log^2 n \big\}$$ we have
$$\begin{aligned}
V_n^2 & =  \sigma^2\, \me \ell_n^{\ssup 2} + O(n^{2/3} \log^3 n).
\end{aligned}$$
Since for $d\ge 4$,
$$\me \ell_n^{\ssup 2} - n\,(2G(0)-1) = O(\log n),$$
we obtain
$$ V_n^2 = n\,\sigma^2\,(2G(0)-1) +  O(n^{2/3} \log^3 n)\, .$$
Thus, if we assume $\sqrt{n} \ll b_n \ll n^{2/3}/\log^{3/2} n=:a_n$,
we have
$$-\frac{b_n^2}{2V_n^2}=- \frac{b_n^2}{2n \sigma^2 (2G(0)-1)} + o(1)\, .$$
Using that
\begin{equation}\label{Gausserr}
1- \Phi(x) = \frac1{\sqrt{2\pi} x}\, e^{-\frac{x^2}{2}}\, \big(1+O(x^{-2})\big), \qquad
\mbox{ as } x\to\infty\, ,
\end{equation}
and abbreviating $\rho_n^2 := 2n \sigma^2 (2G(0)-1)$ we obtain, on the same event,
$$\frac{\frac{V_n}{\sqrt{2\pi} b_n} \, \exp\Big\{ - \frac{b_n^2}{2V_n^2} \Big\}}
{1 - \Phi({b_n}/{\rho_n})} =1+o(1).$$
Therefore, for a constant $c>0$ and all large~$n$,
$$\begin{aligned}
& \left| \frac{{\mathbb P}\{ X_n \ge b_n \} }
{1 - \Phi({b_n}/{\rho_n})} - 1 \right| \\[1mm]
& \le {\mathbb E} \Bigg[ \bigg| \frac{{\mathsf P}\{ \sum \ell_n(z)
\xi(z) \ge b_n \}} { \frac{V_n}{\sqrt{2\pi} b_n} \, \exp\big\{ -
\frac{b_n^2}{2V_n^2} \big\}} - 1 \bigg| \, \, \1\big\{
|\ell_n^{\ssup 2}- \me \ell_n^{\ssup 2}| \le n^{2/3} \log^3 n, \,
\ell_n^{\ssup 3} \le n \log^2 n \big\}\Bigg] +o(1)\\
& \qquad + \prob\big\{ |\ell_n^{\ssup 2}- \me \ell_n^{\ssup 2}| > n^{2/3} \log^3 n \big\}\, e^{c\frac{b_n^2}{n}}
+ \prob\big\{ \ell_n^{\ssup 3} > n \log^2n \big\} \, e^{c\frac{b_n^2}{n}}\, .
\end{aligned}$$
By Proposition~\ref{centredmoment} both probabilities in the last
line are bounded by $\exp\{-cn^{1/3}\}$ if $d\ge 5$, and by 
$\exp\{-c n^{1/3}/\log^{1/2} n \}$ if $d=4$. As $b_n\ll a_n$ we have
$b^2_n/n\ll n^{1/3}$ if $d\ge 5$, and  $b^2_n/n\ll (n/\log^2
n)^{1/3}$ if $d=4$, hence the summands in the last line go to zero,
and together with Proposition~\ref{vitali} this implies
Theorem~\ref{maind>3}. \qed\end{Proof}\medskip

\begin{Proof}{Proof of Proposition~\ref{vitali}.}
Recall Cram\'er's condition~\eqref{cram} and denote $f(h):={\mathbb E} e^{h\xi(0)}$ for all $h\in[0,\theta)$.
For fixed $n\ge 1$ and $h>0$ satisfying the condition
\begin{equation}\label{hcond}
h\, \ell_n^{\ssup{\infty}} \le \sfrac\theta2
\end{equation}
we introduce a family $\{Y_z \colon z\in \Z^d\}$ of independent auxiliary random variables
with distributions
$$\begin{aligned}
P\big\{ Y_z < x \big\} & = \big( f(h \ell_n(z))  \big)^{-1}\,
\int_{-\infty}^x e^{hy}\, d{\mathsf P}\{\ell_n(z)\xi(z)<y\} \, .
\end{aligned}$$
We define
$$\begin{array}{lll}
m_z := E Y_z, & \sigma^2_z := E[(Y_z-m_z)^2], & \gamma_z := E|Y_z-m_z|^3, \\[2mm]
M_n(h) := \sum_{z\in\Z^d} m_z, &
V_n^2(h) := \sum_{z\in\Z^d} \sigma^2_z ,&
\Gamma_n(h) := \sum_{z\in\Z^d} \gamma_z \, . \\[2mm]
\end{array}$$
{F}rom the definition of $Y_z$ we infer that
$${\mathsf P} \big\{ \ell_n(z) \xi(z) < x \big\}
= f(h \ell_n(z))\, \int_{-\infty}^x e^{-hy}\, dP\{Y_z<y\}\, ,$$
and therefore
$${\mathsf P} \Big\{ \sum_{z\in\Z^d} \ell_n(z) \xi(z) \ge b_n \Big\}
= \prod_{z\in\Z^d} f(h \ell_n(z))\, \int_{b_n}^{\infty} e^{-hy}\, dP\Big\{\sum_{z\in\Z^d} Y_z<y\Big\}\, .$$
Substituting $y=M_n(h)+xV_n(h)$ and denoting $T:=(\sum Y_z - M_n(h))/V_n(h)$, we get
\begin{equation}\label{cramrep}\begin{aligned}
{\mathsf P} \Big\{ \sum_{z\in\Z^d} \ell_n(z) \xi(z) \ge b_n \Big\}
=\exp\Big\{ -h M_n(h)& + \sum_{z\in\Z^d} \log f(\ell_n(z)h) \Big\} \\
& \times\,\int_{\frac{b_n-M_n(h)}{V_n(h)}}^{\infty}
\exp\{-hxV_n(h) \} \, dP(T<x).
\end{aligned}\end{equation}
Now we show that \eqref{hcond} implies
that, for some constant~$c>0$, we have
\begin{equation}\label{drei}
h \, V_n^2 - c \, h^3 \,\Gamma_n \le M_n(h) \le h \,  V^2_n
+ c \,h^2\,\Gamma_n\, .
\end{equation}
Obviously,
$$m_z = \frac{\ell_n(z)\, f'(\ell_n(z)h)}{ f(\ell_n(z) h)} \qquad\mbox{ and thus }\qquad
M_n(h) = \sum_{z\in\Z^d} \frac{\ell_n(z)\, f'(\ell_n(z)h)}{ f(\ell_n(z) h)}\, .$$
On the one hand, using that all derivatives of $f$ are increasing, we get
$$f'( \ell_n(z)h ) \le f''(0)\,\ell_n(z) \,h + \sfrac12\,f'''(\ell_n(z)h)\, \ell^2_n(z) \,h^2
\le \sigma^2 \,\ell_n(z)\, h  + \sfrac 12\, f'''(\theta/2) \, \ell_n^2(z) \,h^2,$$
and the second inequality in \eqref{drei} readily follows from this together with the fact that
$f(\ell_n(z)h)\ge 1$.
On the other hand, noting that
$f'(\ell_n(z) h) \ge \sigma^2 \, \ell_n(z) \, h$ and
$$f(\ell_n(z)h ) \le 1 + f'(\ell_n(z) h ) \, \ell_n(z) \, h \le
1+ f'(\theta/2) \, \ell_n(z) \, h \, ,$$
we obtain the bound
$$M_n(h) \ge \sum_{z\in\Z^d} \frac{\sigma^2 \, \ell_n(z) \, h}{1+ f'(\theta/2) \, \ell_n(z) \, h}
= h\,\sigma^2\, \sum_{z\in\Z^d} \ell_n^2(z) - f'(\theta/2)\,\sigma^2\,h^2\, \sum_{z\in\Z^d} \ell_n^3(z)\, .$$
Summarizing, we see that \eqref{drei} holds with $c:=\max\{\sigma^2 f'(\theta/2), \frac 12\,f'''(\theta/2)\}$.

Let $h_n^{\pm}$ denote the positive solutions of the quadratic equations
$$V_n^2\,h \pm c \Gamma_n\, h^2 = b_n\, .$$
It is easy to see that
\begin{equation}\label{vier}
h_n^{\pm} = \frac{b_n}{V_n^2} + O\Big( \frac{\Gamma_n b_n^2}{V_n^6} \Big) \qquad\mbox{ as } n\to\infty\, ,
\end{equation}
provided that $\Gamma_nb_n= O(V_n^4)$.

{F}rom our assumption $\Gamma_n \le n \log^2 n$ we get
$\ell_n^{(\infty)} \le n^{1/3} \log^{2/3} n$\and thus \eqref{hcond}
holds for all $h\le \theta/ (2n^{1/3} \log^{2/3} n)$. Since $b_n \le
n^{2/3}/ \log n$ and $\Gamma_n b_n^2 \le n^{7/3}$ but $V_n^2 \ge n$
we obtain that $h_n^- \le n^{-1/3}/ \log n + O(n^{-2/3})$ and thus
$h_n^-$ is in the domain given by \eqref{hcond}, for all large~$n$.
Hence the inequalities~\eqref{drei} hold for all $0<h\le h_n^-$ and
so, on the one hand, we have $M(h_n^-)\ge b_n$, and on the other
hand, as $h_n^+<h_n^-$, we have $M(h_n^+)\le b_n$. Therefore there
exists $h_n\in [h_n^+, h_n^-]$ such that $M(h_n)=b_n$. Applying
\eqref{vier} gives
\begin{equation}\label{fuenf}
h_n= \frac{b_n}{V_n^2} + O\Big( \frac{\Gamma_n b_n^2}{V_n^6} \Big) \qquad\mbox{ as } n\to\infty\, .
\end{equation}
Clearly,
$$\log f\big( \ell_n(z)h_n \big) = \log\big( 1 + \sfrac{\sigma^2}2\,\ell^2_n(z)\,h_n^2 +
O(\ell_n^3(x)h_n^3) \big) = \sfrac{\sigma^2}2\,\ell^2_n(z)\,h_n^2 + O(\ell_n^3(x)h_n^3)\, .$$
Thus, in view of \eqref{fuenf},
\begin{equation}\label{sechs}
-h_n M_n(h_n) + \sum_{z\in\Z^d} \log f\big( \ell_n(z) h_n\big)
= - h_n\,b_n + \sfrac 12\, V_n^2\,h_n^2 + O\big( \Gamma_n h_n^3 \big) = - \frac{b_n^2}{2 V_n^2} +
O \Big( \frac{\Gamma_n b_n^3}{V_n^6} \Big) \, .
\end{equation}
Putting $h=h_n$ in \eqref{cramrep} and using \eqref{sechs}, we obtain
\begin{equation}\label{sieben}
 {\mathsf P} \big\{ \ell_n(z) \xi(z) \ge b_n \big\}
=\exp\Big\{ - \frac{b_n^2}{2V_n} + O\Big(\frac{\Gamma_n b_n^3}{V_n^6}\Big)\Big\}
\,\int_{0}^{\infty} e^{-x h_n V_n(h_n)}  \, dP(T<x).
\end{equation}
Integrating by parts gives, for a standard normal random variable~$N$,
$$\begin{aligned}
\int_0^\infty & e^{-x h_n V_n(h_n)} \, dP\{T<x\}
= \int_0^\infty P\{T<x\}\, h_n \, V_n(h_n)\, e^{-h_n V_n(h_n)\,x}\, dx\\
& = \int_0^\infty P\{N<x\} \, h_n \, V_n(h_n) \, e^{- h_n\, V_n(h_n)\,x}\, dx
+\int_0^\infty \Delta(x) \, h_n\, V_n(h_n) e^{-h_n\, V_n(h_n)\, x}\, dx \\
& = \frac1{\sqrt{2\pi}} \int_0^\infty \exp\big\{ - h_n \,V_n(h_n) \, x - \sfrac{x^2}{2} \big\} \, dx
+ \int_0^\infty \Delta(x)\, h_n\,V_n(h_n)\, e^{-h \, V_n(h_n)\, x}\, dx,
\end{aligned}$$
where $\Delta(x):=P\{T<x\} - P\{N<x\}$. By Esseen's inequality, see for example~\cite[Theorem V.3]{Pe75},
there exists an abolute constant~$C>0$, such that
$$\sup_x |\Delta(x)| \le C \, \frac{\Gamma_n(h_n)}{V_n^3(h_n)}\, .$$
Therefore
$$\Big|  \int_0^\infty e^{-x h_n V_n(h_n)} \, dP\{T<x\}
-  \frac1{\sqrt{2\pi}} \int_0^\infty \exp\big\{ - h_n \,V_n(h_n) \, x - \sfrac{x^2}{2} \big\} \, dx \Big|
\le C\,\frac{\Gamma_n(h_n)}{V_n^3(h_n)}\, .$$
Evidently,
\begin{equation}\label{acht}\begin{aligned}
\frac1{\sqrt{2\pi}}   \int_0^\infty & \exp\big\{ - h_n \,V_n(h_n) \, x - \sfrac{x^2}{2} \big\} \, dx \\
& = \frac{1}{\sqrt{2\pi}}\, \exp\Big\{ \frac{h_n^2 V_n^2(h_n)}{2} \Big\} \,
\int_0^\infty \exp\Big\{ -\frac{(x+h_n V_n(h_n))^2}{2} \Big\} \, dx \\
& = \exp\Big\{ \frac{h_n^2 V_n^2(h_n)}{2} \Big\}\, \Big( 1- \Phi\big(h_nV_n(h_n)\big) \Big)\, .
\end{aligned}\end{equation}
We now show that, for a suitable constant~$C>0$,
\begin{equation}\label{perturbmom}
V_n^2(h_n) = V_n^2 + O(\Gamma_n h_n)\qquad \mbox{ and } \qquad\Gamma_n(h_n) \le C\, \Gamma_n.
\end{equation}
\emph{First},  we obtain that
$$\begin{aligned}
V_n^2(h_n) & = \sum_{z\in\Z^d} \sigma^2_z =  \sum_{z\in\Z^d}
\ell^2_n(z)\, \frac{f''(\ell_n(z)h_n) - (f'(\ell_n(z) h_n))^2 }{f(\ell_n(z)h_n)} \\
& =  \sum_{z\in\Z^d} \ell^2_n(z)\, \big(\sigma^2 + O(\ell_n(z)\, h_n) \big) = V_n^2 + O(\Gamma_n h_n)\, .
\end{aligned}$$
\emph{Second}, for an upper estimate of $\Gamma_n(h_n)$, we note that
$$E|Y_z|^3 = 2\int_{-\infty}^0 |y|^3\, dP\{Y_z<y\} + E Y_z^3\, .$$
{F}rom the definition of $Y_z$ we get, on the one hand,
$$\begin{aligned}
\int_{-\infty}^0|y|^3\, dP\{Y_z<y\} &
= \sfrac{1}{f(\ell_n(z)h)}\, \int_{-\infty}^0 |y|^3 \,e^{hy} \, d{\mathsf P}\big\{ \xi(z) < \sfrac y{\ell_n(z)}\}\\
& \le \ell_n^3(z)\, \int_{-\infty}^0 |x|^3 \, d{\mathsf P}\big\{ \xi(z) < x\}
\le \ell_n^3(z)\, {\mathsf E}|\xi(0)|^3\, ,
\end{aligned}$$
and, on the other hand,
$$E Y_z^3 = \frac{f'''(\ell_n(z)h) \, \ell_n^3(z)}{f(\ell_n(z) h)}
\le \ell_n^3(z)\, f'''(\theta/2)\, .$$
The two bounds imply that
$E|Y_z|^3 \le \big( f'''(\theta/2) + 2 \gamma \big) \ell_n^3(z),$
and combining this with $m_z\le f'(\theta/2)\, \ell_n(z)$ gives
$\gamma_z\le E|Y_z|^3 + m_z^3 \le C \ell^3_n(z)$ and therefore we
have proved~\eqref{perturbmom}.

{F}rom \eqref{fuenf} and \eqref{perturbmom} we thus get
$$\begin{aligned}
h_n V_n(h_n) & = \Big( \frac{b_n}{V_n^2} + O\big( \sfrac{\Gamma_n b_n^2}{V_n^6}\big)\Big)\,
\Big( V_n^2 + O\big(\sfrac{\Gamma_nb_n}{V_n^2} \big)\Big)^{1/2}
 = \frac{b_n}{V_n}\, \Big( 1+ O \big( \sfrac{\Gamma_n b_n}{V_n^4} \big) \Big)\, .
\end{aligned}$$
Recalling that $b_n \gg \sqrt{n}$ and $V_n^2\le An$ we conclude that
$h_n V_n(h_n) \to \infty$. Then, using \eqref{Gausserr},
$$\begin{aligned}
e^{h_n^2 V_n^2(h_n)/2} \big( 1- \Phi(h_n V_n(h_n) \big)
& = \frac{1}{\sqrt{2\pi} h_n V_n(h_n)} \, \Big(1+ O\big( \sfrac1{h_n^2 V_n^2(h_n)} \big) \Big) \\
& = \frac{V_n}{\sqrt{2\pi} b_n} \, \Big( 1 + O\big( \sfrac{\Gamma_n b_n}{V_n^4} \big)
+ O\big( \sfrac{V_n^2}{b_n^2} \big) \Big)\, .
\end{aligned}$$
Substituting this into \eqref{acht} gives
$$\int_0^\infty e^{-h_n V_n(h_n) x}\, dP\{T<x\}
= \frac{V_n}{\sqrt{2\pi} b_n} \, \Big( 1 + O\big( \sfrac{\Gamma_n b_n}{V_n^4} \big)
+ O\big( \sfrac{V_n^2}{b_n^2} \big) \Big)\, ,$$
and the result follows by plugging this into \eqref{sieben}.
\qed\end{Proof}

\section{Moderate deviations in dimensions $d\ge 3$: Proof of Theorem~\ref{maind>2}}

\subsection{Proof of the upper bound in Theorem~\ref{maind>2}}

We fix $\epsilon>0$ and let $A:= 2G(0)-1 + 3\epsilon.$
Our aim is to show that
\begin{equation}\label{uppertarget}
\limsup_{n\to\infty} \frac{n}{b_n^2} \log \prob\Big\{ \sum_{z\in \Z^d} \ell_n(z) \xi(z) \ge b_n \Big\} \le
-\frac{1}{ 2\sigma^2\,A } .
\end{equation}
We note that, for any fixed $\eta>0$,
\begin{equation}\begin{aligned}\label{splitup}
\prob\Big\{ \sum_{z\in \Z^d} \ell_n(z) \xi(z) \ge b_n  \Big\}
& \le  \prob\Big\{ \sum_{z\in \Z^d} \ell_n(z) \xi(z)\ge b_n, \,\,
\ell_n^{\ssup \infty} \le \sfrac{\eta n}{b_n},  \, \ell_n^{\ssup 2} \le A\,n  \Big\}\\
& \qquad \qquad\qquad + \prob\Big\{ \ell_n^{\ssup \infty} \ge \sfrac{\eta n}{b_n} \Big\}
+ \prob\big\{\ell_n^{\ssup 2} \ge A\, n \big\}.
\end{aligned}\end{equation}
To see that the second summand is negligible apply
Lemma~\ref{maxdev} with $a_n=\eta n/b_n$, which gives
\begin{equation}\label{second}
\limsup_{n\to\infty}\frac{n}{b_n^2}\log \prob\big\{ \ell_n^{\ssup \infty} > \sfrac{\eta n}{b_n} \big\} \le
\limsup_{n\to\infty}\frac{n \log n}{b_n^2} -c \frac{\eta n^2}{b_n^3} = - \infty.
\end{equation}
To see that the third term in \eqref{splitup} is negligible, recall from~\eqref{expd} that
$\me \ell_n^{\ssup 2}\sim n (2G(0)-1)$ and therefore, for all large $n$,
$$\begin{aligned}
\prob\big\{ \ell^{\ssup 2}_n \ge A\, n \big\} &
\le \prob \Big\{ \ell^{\ssup 2}_n-\me \ell^{\ssup 2}_n
\ge \big( \sfrac{A-1-\epsilon}{2}-G(0)\big) n \Big\}
 = \prob \big\{ \ell^{\ssup 2}_n-\me \ell^{\ssup 2}_n  \ge \epsilon n \big\}.
\end{aligned}$$
{F}rom Proposition~\ref{centredmoment} we know that for $b_n\ll n^{2/3}$, if $d\ge 4$,
$$\limsup_{n\to\infty}\frac{n}{b_n^2}\,\log \prob \big\{ \ell_n^{\ssup 2}-\me \ell_n^{\ssup 2}  \ge \epsilon n \big\}
\le \limsup_{n\to\infty} -c \frac{n^{3/2}}{b_n^2 \log n} = -\infty,$$
and, if $d= 3$,
$$\limsup_{n\to\infty}\frac{n}{b_n^2}\,\log \prob \big\{ \ell_n^{\ssup 2}-\me \ell_n^{\ssup 2}  \ge \epsilon n \big\}
\le \limsup_{n\to\infty} -c \frac{n^{4/3}}{b_n^2} = -\infty.$$
Combining this, we get
\begin{equation}\begin{aligned}\label{third}
\limsup_{n\to\infty}\frac{n}{b_n^2}\log \prob\big\{ \ell^{\ssup 2}_n \ge A\, n \big\}
 = - \infty.
\end{aligned}\end{equation}
It remains to investigate the first term on the right hand side of \eqref{splitup}.
For this purpose, for the moment fix $\{ \ell_n(z) \, : \, z\in\Z^d\}$ such that
$$\ell_n^{\ssup \infty} \le \frac{\eta n}{b_n}\qquad
\mbox{ and } \qquad \ell^{\ssup 2}_n \le A n,$$ and just look at probabilities
for the i.i.d. variables $\{\xi(z) \, : \,z\in \Z^d\}$. Denote
$f(h):= {\sf E} e^{h \xi(0)}$ for all $h<\theta$,
which is well-defined by Cram\'er's condition. Recall that
$$f(h) = \exp\big\{ \sfrac{1}{2} \, h^2\, \sigma^2 (1+o(h)) \big\}  \qquad\mbox{ as } h\downarrow 0.$$
In particular, given any $\delta>0$, we may choose a small $\eta>0$ such that
\begin{equation}\label{useful}
f\Big( \frac{ b_n\, \ell_n(x)}{ \sigma^{2}\, \ell^{\ssup 2}_n} \Big)
\le \exp\Big\{ (1+\delta) \, \frac{b_n^2\ell_n^2(x)}{2\sigma^{2}\,(\ell_n^{\ssup 2})^2 }\Big\} ,
\end{equation}
where we use that $b_n \ell_n(x)/\ell^{\ssup 2}_n\le \eta$.
{F}rom Chebyshev's inequality and independence we get that
$$\begin{aligned}
{\sf P}\Big\{ &  \sum_{x\in \Z^d} \ell_n(x) \xi(x) \ge b_n \Big\} \le
\prod_{x\in\Z^d} f\Big( \frac{ b_n \ell_n(x)}{ \sigma^{2}\, \ell^{\ssup 2}_n} \Big) \,
\exp\Big\{ - \frac{b_n^2}{ \sigma^{2}\, \ell^{\ssup 2}_n}\Big\} \\
& \le \exp\Big\{ (1+\delta) \, \frac{b_n^2}{2\sigma^{2}\,\ell^{\ssup 2}_n }\Big\} \,
\exp\Big\{ - \frac{b_n^2}{ \sigma^{2}\, \ell^{\ssup 2}_n}\Big\}
 \le \exp\Big\{ -\big(1-\delta\big) \, \frac{b_n^2}{2\sigma^{2}\, A\, n }\Big\} .
\end{aligned}$$
We can now average over the random walk again, and get \eqref{uppertarget} {f}rom \eqref{splitup} together
with \eqref{second} and \eqref{third}, recalling that $\delta>0$ was arbitrary. This completes the proof.\qed
\bigskip

\subsection{Proof of the lower bound in Theorem~\ref{maind>2}}\label{LBd>2}

We impose `typical behaviour' on $\ell^{\ssup 2}_n$ and $\ell^{\ssup \infty}_n$.
More precisely, fix an arbitrary $\epsilon\in(0,1)$, and also fix
$\eta>0$ which we specify later. We have
\begin{equation}\label{split}\begin{aligned}
\prob\Big\{ \sum_{z\in \Z^d} \ell_n(z) \xi(z) \ge b_n  \Big\} &
\ge  \prob\Big\{ \sum_{z\in \Z^d} \ell_n(z) \xi(z) \ge b_n, \,\, \ell^{\ssup \infty}_n \le \sfrac{\eta n}{b_n},
\,\ell^{\ssup 2}_n \le A\,n  \Big\}\\
& = \me\Big\{ {\sf P}\Big\{ \sum_{z\in\Z^d} \ell_n(z)\, \xi(z) \ge b_n\Big\}\,
\1\{ \ell_n^{\ssup \infty} \le \sfrac{\eta n}{b_n}, \,\ell^{\ssup 2}_n \le A\,n \} \Big\},
\end{aligned}\end{equation}
where $A:= 2G(0)-1 + 3\epsilon$ and ${\sf P}$ refers to the probability with respect to the
scenery only. To study the inner probability we now suppose that, for the moment, a random
walk sample is fixed, such that
$$\ell_n^{\ssup \infty} \le \sfrac{\eta n}{b_n}\qquad
\mbox{ and } \qquad \ell^{\ssup 2}_n \le A n.$$
Denote $\gamma:=\me|\xi(0)|^3<\infty$. Hence
the variance of the random variable $\sum_{z\in \Z^d} \ell_n(z) \xi(z)$ with respect to $\sf P$
is given by $V_n^2:= \sigma^2 \sum_{z\in\Z^d} \ell_n^2(z)$ and the Lyapunov ratio by
$L_n := \gamma\,V_n^{-3} \sum_{z\in \Z^d} \ell_n^3(z)$.
By \cite[Theorem~2]{Na02} there exist constants $c_1, c_2>0$ such that, for all
$\frac 32 V_n\le x\le  \frac{V_n}{196 L_n}$,
\begin{equation}\label{4.1.'}
{\sf P}\Big\{ \sum_{z\in \Z^d} \ell_n(z) \xi(z) \ge x \Big\}
\ge \big(1-\Phi(\sfrac{x}{V_n})\big)\, \exp\big\{ -c_1 x^3\, L_n V_n^{-3} \big\}
\, \big( 1- c_2 x L_n V_n^{-1} \big) .
\end{equation}
Now suppose that $\eta>0$ is chosen to satisfy the three inequalities
 $$\eta<\sigma^4/(196\gamma), \, c_1\eta\gamma \sigma^{-6}<\epsilon, \, \mbox{ and }
c_2\eta\gamma \sigma^{-4} < \epsilon.$$
Using the upper bound on $\ell_n^{\ssup \infty}$, we get that
$L_n \le \frac{\gamma\eta n}{\sigma^2b_n} V_n^{-1}$. Therefore,
$${\sf P}\Big\{ \sum_{z\in \Z^d} \ell_n(z) \xi(z) \ge x \Big\}
\ge \big(1-\Phi(\sfrac{x}{V_n})\big)\, \exp\big\{ -c_1 \,\eta\, \sfrac{x^3}{b_n\, \sigma^2} \, n \,V_n^{-4} \big\}
\, \big( 1- c_2\, \gamma \,\eta \, \sfrac{x}{b_n\, \sigma^2}\, n \,V_n^{-2} \big) ,$$
for all $(3/2) V_n\le x\le  (b_n V_n)/(196 \eta n)$.
We can use this inequality for $x=b_n$. Indeed, as $V_n^2\le A\sigma^2\,n$ we get $b_n\ge (3/2) V_n$,
if $n$ exceeds some constant depending only on $\sigma^2$.
Also $V_n^2\ge \sigma^2 n$ and $\eta<\sigma^4/(196\gamma)$, therefore
$$b_n\le b_n \sigma^2 V_n^2 / (196\gamma \eta n) \le V_n/(196 L_n).$$
Hence,
\begin{equation}\label{mainstep}
{\sf P}\Big\{ \sum_{z\in \Z^d} \ell_n(z) \xi(z) \ge b_n \Big\}
\ge \big(1-\Phi(\sfrac{b_n}{V_n})\big)\, \exp\big\{ -c_1\, \eta \, \gamma\, \sigma^{-6}\,\sfrac{b_n^2}{n} \big\}
\, \big( 1- c_2 \gamma \, \sigma^{-4}\, \eta  \big) .
\end{equation}
Substituting \eqref{mainstep} into \eqref{split} gives
\begin{equation}\label{L4}\begin{aligned}
\prob\Big\{ & \sum_{z\in \Z^d} \ell_n(z) \xi(z) \ge b_n  \Big\}\\
& \ge  \big( 1- c_2 \gamma \sigma^{-4}\,\eta  \big) \,\exp\big\{ -c_1\, \gamma\sigma^{-6}\,
\eta \,\sfrac{b_n^2}{n} \big\} \,
\me\Big[ \big(1-\Phi(\sfrac{b_n}{V_n})\big) \, \1\big\{ V_n^2 \le A\sigma^2\,n, \, \ell_n^{\ssup \infty}
\le \sfrac{\eta n}{b_n} \big\} \Big] \\
& \ge \big( 1- \epsilon \big) \,\exp\big\{ -\epsilon\,\sfrac{b_n^2}{n} \big\} \,
\me\Big[ \big(1-\Phi(\sfrac{b_n}{V_n})\big) \, \1\big\{ V_n^2 \le A\sigma^2\,n \big\} \Big]-
\prob\Big\{ \ell_n^{\ssup \infty} \ge \sfrac{\eta n}{b_n} \Big\} .
\end{aligned}\end{equation}
Since, by a standard estimate,
$(1-\Phi(z)\big) \ge \exp\{ -(1+\eta)\,z^2/2\}$ for all sufficiently large~$z$, we get
\begin{equation}\label{L5}\begin{aligned}
\me\Big[ \big(1-\Phi(\sfrac{b_n}{V_n})\big) \, \1\big\{ V_n^2 \le A\sigma^2\,n \big\} \Big]
& \ge \me \exp\Big\{ -\frac{(1+\eta) b_n^2}{2V_n^2} \Big\}  - \prob\big\{ V_n^2 \ge A\sigma^2\,n \big\} .
\end{aligned}\end{equation}
By Jensen's inequality, we obtain
$$ \me \exp\Big\{ -\frac{(1+\eta) b_n^2}{2V_n^2} \Big\}
\ge \exp\Big\{ -\frac{(1+\eta) b_n^2}{2\sigma^2n}\,\me \frac{n}{\ell^{\ssup 2}_n}  \Big\}.$$
Using Proposition~\ref{centredmoment} and the Borel-Cantelli lemma,
$$\lim_{n\to\infty} \frac {\ell^{\ssup 2}_n}n  = \lim_{n\to\infty} \frac {\me \ell^{\ssup 2}_n}n
= 2G(0)-1 \quad \mbox{ almost surely,}$$
and using further that $n/\ell^{\ssup 2}_n\le 1$, we obtain that
$$\lim_{n\to\infty} \me \frac{n}{\ell^{\ssup 2}_n} = \frac 1{2G(0)-1}.$$
Then, for all $n$ sufficiently large,
\begin{equation}\label{L6}\begin{aligned}
\me\Big[ \big(1-\Phi(\sfrac{b_n}{V_n})\big)\Big]
& \ge \exp\Big\{ -\frac{(1+2\eta) b_n^2}{2\sigma^2n (2G(0)-1-\epsilon)}\Big\} .
\end{aligned}\end{equation}
Combining \eqref{L4}, \eqref{L5} and \eqref{L6} gives
$$\begin{aligned}
\prob\Big\{ \sum_{z\in \Z^d} \ell_n(z) \xi(z) \ge b_n  \Big\}
\ge (1-\epsilon)\, & \, \exp\Big\{ -\Big( \epsilon+\frac{1+2\eta}{2\sigma^2(2G(0)-1-\epsilon)}\Big) \frac{b_n^2}{n}\Big\}\\
& \qquad - \prob\Big\{ \ell_n^{\ssup \infty} \ge \sfrac{\eta n}{b_n} \Big\}  -
\prob\big\{ \ell^{\ssup 2}_n  \ge A\, n \big\} .
\end{aligned}$$
The required lower bound follows {f}rom the estimates \eqref{second} and \eqref{third} for the subtracted probabilities,
and the fact that $\epsilon>0$ can be chosen arbitrarily small, whence $\eta$ also becomes arbitrarily small.\qed

\section{Moderate deviations in dimension $d=2$: Proof of Theorem~\ref{maind=2}}

We use the following moderate deviation principle
for the self-intersection local time in the planar case, which is due to Bass,
Chen and Rosen~\cite[Theorem~1.1 and (3.2)]{BCR06}: If $x_n\to\infty$ and $x_n=o(n)$,
then for every $\lambda>0$,
\begin{equation}
\label{p1}
\lim_{n\to\infty}\frac{1}{x_n}\log\mathbb{P}\big\{ \ell_n^{\ssup 2}-\mathbb{E}\ell_n^{\ssup 2}
\geq \lambda\, nx_n\big\}= \lim_{n\to\infty}\frac{1}{x_n}\log\mathbb{P}\big\{
|\ell_n^{\ssup 2}-\mathbb{E}\ell_n^{\ssup 2} | \geq \lambda\, nx_n\big\}
= - \frac{\lambda\sqrt{\det \Gamma}}{2\varkappa^4},
\end{equation}
where again $\varkappa$ is the optimal constant in the Gagliardo-Nirenberg inequality.

\subsection{Proof of Theorem~\ref{maind=2}(a)}

The proof is largely analogous to that of Theorem~\ref{maind>2} replacing Proposition~\ref{centredmoment}
by \eqref{p1}.
Starting with the \emph{upper bound}, for any fixed $\epsilon>0$, we use the decomposition
$$\begin{aligned}
\prob\Big\{ \sum_{z\in \Z^2} \ell_n(z) \xi(z) \ge b_n  \Big\}
& \le  \prob\Big\{ \sum_{z\in \Z^2} \ell_n(z) \xi(z)\ge b_n, \,\,
\ell_n^{\ssup \infty} \le \sfrac{\sqrt{n}(\log n)^5}{b_n}, \, \ell_n^{\ssup 2} \le A\,n \log n \Big\}\\
& \qquad \qquad\qquad + \prob\Big\{ \ell_n^{\ssup \infty} \ge \sfrac{\sqrt{n}(\log n)^5}{b_n} \Big\}
+ \prob\big\{ \ell_n^{\ssup 2} \ge A\, n \log n\big\},
\end{aligned}$$
where $A:=(\pi\sqrt{\det \Gamma})^{-1}+4\epsilon$. The estimate for the last probability follows
{f}rom \eqref{p1}. Indeed, by~\eqref{expd}, for sufficiently large~$n$,
$$\begin{aligned}
\prob\big\{ \ell_n^{\ssup 2} \ge A\, n \log n\big\} &
\le \prob\big\{ \ell_n^{\ssup 2} - \me \ell_n^{\ssup 2} \ge
\big(A-(\pi\sqrt{\det \Gamma})^{-1}-\epsilon\big)\, n \log n\big\}
\le  n^{-\epsilon\,\sqrt{\det \Gamma}\,\varkappa^{-4}},
\end{aligned}$$
hence, as $b_n\ll n^{\frac 12} \log n$,
\begin{equation}\label{midval}
\limsup_{n\to\infty}\frac{n \log n}{b_n^2} \log \prob\big\{ \ell_n^{\ssup 2} \ge A\, n \log n\big\}
= - \infty.\end{equation}
Moreover, applying Lemma~\ref{maxdev}, we get
\begin{equation}\label{maxival}\begin{aligned}
\limsup_{n\to\infty}\frac{n \log n}{b_n^2}
\log & \, \prob\big\{ \ell_n^{\ssup \infty} > b_n^{-1}\,\sqrt{n} \,(\log n)^5 \big\} \\
& \le
\limsup_{n\to\infty}\frac{n (\log n)^2}{b_n^2} -c \frac{n^{\frac 32} (\log n)^4}{b_n^3} = - \infty.
\end{aligned}\end{equation}
We now look at fixed local times $\{\ell_n(z) \colon z\in\Z^2\}$ satisfying
the conditions $\max \ell_n(z) \le b_n^{-1}\,\sqrt{n}(\log n)^5$ and
$\ell_n^{\ssup 2} \le A\,n \log n$. Note that, together with the trivial
inequality $\ell_n^{\ssup 2}\ge n$, this implies
$$\lim_{n\uparrow\infty} \frac{b_n \ell_n(z)}{\sigma^2\,\ell_n^{\ssup 2}} = 0 .$$
Hence, for arbitrary $\delta>0$, if $n$ is sufficiently large, an application of Chebyshev's inequality and the
estimate \eqref{useful} for the Laplace transform $f$ of $\xi(z)$, gives, for $n$ larger than some absolute constant,
$$\begin{aligned}
{\sf P}\Big\{  \sum_{z\in \Z^2} \ell_n(z) \xi(z) \ge b_n \Big\} & \le
\prod_{z\in\Z^d} f\big( \sfrac{ b_n \ell_n(z)}{ \sigma^{2}\, \ell^{\ssup 2}_n} \big) \,
\exp\big\{ - \sfrac{b_n^2}{ \sigma^{2}  \, \ell^{\ssup 2}_n}\big\}
 \le \exp\big\{ -(1-\delta) \, \sfrac{b_n^2}{2\sigma^{2}\, A\, n\log n }\big\} .
\end{aligned}$$
Averaging over the local times again, we obtain
$$\begin{aligned}
\limsup_{n\uparrow\infty}\sfrac{n \log n}{b_n^2} \log \prob\Big\{ \sum_{z\in \Z^2} \ell_n(z) \xi(z)\ge b_n, \,
\ell_n^{\ssup \infty} \le \sfrac{\sqrt{n}(\log n)^5}{b_n}, \, \ell_n^{\ssup 2} \le A\,n \log n \Big\}
& \le \sfrac{-(1-\delta)}{2\sigma^{2}\, A},
\end{aligned}$$
so that the claimed upper bound follows, as $\epsilon,\delta>0$ were arbitrary.\medskip

Turning to the \emph{lower bound}, we fix $\epsilon>0$ again, and use that
\begin{equation}\label{spl}\begin{aligned}
\prob\Big\{ \sum_{z\in \Z^2} \ell_n(z) \xi(z) \ge b_n  \Big\}
\ge  \me\Big\{ {\sf P}\big\{ & \, \sum_{z\in\Z^2} \ell_n(z)\, \xi(z) \ge b_n\big\}\, \\
& \quad \times \1\big\{ \ell_n^{\ssup \infty} \le \sfrac{\sqrt{n}(\log n)^5}{b_n},
 \, \ell_n^{\ssup 2} \le A\,n \log n\big\} \Big\},
\end{aligned}\end{equation}
where $A:=(\pi\sqrt{\det \Gamma})^{-1}+4\epsilon$.
To obtain a lower bound for the inner probability we argue as in Theorem~\ref{maind>2},
relying on the estimates of \cite[Theorem 2]{Na02}. This gives
$$\begin{aligned}
{\sf P}\Big\{ & \, \sum_{z\in \Z^d} \ell_n(z) \xi(z) \ge b_n \Big\}
 \ge \big(1-\Phi(\sfrac{b_n}{V_n})\big)\, \exp\big\{ -c_1\, \gamma\sigma^{-6}\, b_n^2 \, n^{-\frac 32}\,
(\log n)^3 \big\} \, \big( 1- c_2 \, \gamma\sigma^{-4}\, n^{-\frac12} (\log n)^4  \big) .
\end{aligned}$$
We now show that
\begin{equation}\label{expelim}
\lim_{n\uparrow\infty} \me\Big[\frac{n \log n}{\ell_n^{\ssup 2}}\Big]
= \pi\,\sqrt{\det\Gamma} \,.
\end{equation}
For this purpose define the random variables
$Y_n := \frac 1{n}\, \ell_n^{\ssup 2} - (\pi\sqrt{\det\Gamma})^{-1} \,\log n$
and note that
$$\frac{n \log n}{\ell_n^{\ssup 2}} = \pi\sqrt{\det\Gamma} -
\pi\sqrt{\det\Gamma} \, \frac{ Y_n}{\frac 1{n}\, \ell_n^{\ssup
2}}.$$ It suffices to show that the expectation of the fraction on
the right converges to zero. As $|Y_n| \le \eps \log n$ implies that
$\sfrac 1n\ell_n^{\ssup 2} \ge ((\pi\sqrt{\det\Gamma})^{-1}-\eps)
\log n$ we obtain, for any small $\eps>0$, that
\begin{equation}\label{1teil}
\me \Big[ \frac{|Y_n|}{\frac 1{n}\, \ell_n^{\ssup 2}} \, \1\{ |Y_n| \le \eps \log n \} \Big]
\le \frac{\eps}{(\pi\sqrt{\det\Gamma})^{-1} - \eps}\, .
\end{equation}
Also, as $\sfrac 1n\ell_n^{\ssup 2}\ge 1$ and using \eqref{p1} with $\lambda=\eps$ and $x_n=\log n$, for any $0<\eps<\delta$,
\begin{equation}\label{2teil}
\me \Big[ \frac{|Y_n|}{\frac 1{n}\, \ell_n^{\ssup 2}}\, \1\{ \eps \log n < |Y_n| \le \delta \log n \} \Big]
\le \delta\, (\log n) \, \prob\{ |Y_n| > \eps\, \log n \} \longrightarrow 0 ,
\end{equation}
and, using \eqref{p1} with $\lambda=\delta$ and $x_n=\log n$, if $\delta>0$ is sufficiently large,
\begin{equation}\label{3teil}
\me \Big[ \frac{|Y_n|}{\frac 1{n}\, \ell_n^{\ssup 2}} \, \1\{ |Y_n| > \delta \log n \} \Big]
\le  n \,  \prob\{ |Y_n| > \delta\, \log n \} \longrightarrow 0 \, .
\end{equation}
We obtain that $\lim \me|Y_n|/\frac 1n \ell_n^{\ssup 2}=0$, and hence \eqref{expelim}, by combining \eqref{1teil}, \eqref{2teil}, and \eqref{3teil}.

Repeating the arguments of the $d\ge 3$ case, given in Section~\ref{LBd>2}, gives
$$\begin{aligned}
\prob\Big\{ & \sum_{z\in \Z^d} \ell_n(z) \xi(z) \ge b_n  \Big\} \\
& \ge \big( 1- c_2 \, \gamma\sigma^{-4} \, n^{-\frac12} (\log n)^4  \big) \,
\exp\big\{ -c_1\, \gamma\sigma^{-6}\, b_n^2 \, n^{-\frac 32}\, (\log n)^3 \big\} \,
\exp\Big\{-\frac{(1+\eps)^2 \pi b_n^2}{2\sigma^2 n \log n}\Big\} \\
& \qquad\quad- \prob\big\{ \ell_n^{\ssup \infty} \ge \sfrac{\sqrt{n}(\log n)^5}{b_n} \big\}  -
\prob\big\{ \ell_n^{\ssup 2}  \ge A\, n\log n \big\} .
\end{aligned}$$
The result follows, by observing that the first two factors on the right
converge to one, recalling \eqref{maxival}, \eqref{midval} and that $\epsilon>0$ was arbitrary.\qed

\subsection{Proof of Theorem~\ref{maind=2}(b)}

Again, we start with the \emph{upper bound}.
Since $\mathbb{E} \ell_n^{\ssup 2} \sim (\pi\sqrt{\det\Gamma})^{-1}n\log n$, we can conclude {f}rom~(\ref{p1})
that, for $\log n\ll x_n\ll n$,
\begin{equation}
\label{p2}
\lim_{n\to\infty}\frac{1}{x_n}\log\mathbb{P}\{\ell_n^{\ssup 2}\geq\lambda nx_n\}=
-\frac{\lambda}{2\varkappa^{4}}\,\sqrt{\det\Gamma}\, .
\end{equation}

For arbitrary $N\geq1$ and $0<\delta<1$,
\begin{eqnarray}
\label{p3}
\mathbb{P}\{X_n\geq b_n\}\leq
\sum_{i=0}^{N-1}\mathbb{P}\bigl\{X_n\geq b_n, \,  \ell_n^{\ssup 2} \in(i\delta a_n,(i+1)\delta a_n]\bigr\}
+\mathbb{P}\{\ell_n^{\ssup 2} >N\delta a_n\},
\end{eqnarray}
where $a_n:=b_n\sqrt{n}$. Note that $a_n\gg n\log n$. Hence, in view of (\ref{p2}),
\begin{equation}
\label{p4}
\mathbb{P}\{ \ell_n^{\ssup 2} >N\delta a_n\}\leq\exp\Bigl\{-\frac{N\delta a_n \sqrt{\det\Gamma}}
{3\varkappa^{4}n}\Bigr\}
\end{equation}
for all sufficiently large $n$. Fix $i\geq1$ and $\eta\in (0, \theta\sigma^2)$. Then,
$$\begin{aligned}
\mathbb{P} & \,\bigl\{X_n\geq b_n, \, \ell_n^{\ssup 2}\in(i\delta a_n,(i+1)\delta a_n]\bigr\} \\
& \le \mathbb{P}\Bigl\{X_n\geq b_n, \,  \ell_n^{\ssup 2}\in(i\delta a_n,(i+1)\delta a_n], \,
\ell_n^{\ssup \infty} \le \eta i \delta \sqrt{n} \Bigr\}
+ \mathbb{P} \big\{ \ell_n^{\ssup \infty} > \eta i \delta \sqrt{n} \big\}\, .
\end{aligned}$$
Using Lemma~\ref{maxdev}, we get
\begin{equation}
\label{p5}
\mathbb{P}\big\{\ell_n^{\ssup \infty}>\eta i \delta \sqrt{n}\big\}
\leq\exp\Bigl\{-c\frac{\eta i \delta \sqrt{n}}{\log n}\Bigr\}.
\end{equation}
On the event $\big\{\ell_n^{\ssup \infty}\leq\eta i \delta \sqrt{n}, \ell_n^{\ssup 2}\in(i\delta a_n,(i+1)\delta a_n]\big\}$,
we obtain,
$$\frac{b_n \ell_n(z)}{\sigma^2 \ell_n^{\ssup 2}} \le
\frac{b_n \eta i \delta \sqrt{n}}{\sigma^2 i \delta a_n} =
\frac{\eta}{\sigma^2} < \theta. $$
Therefore, we can use Chebyshev's inequality as before, which gives
$$
\mathsf{P}\Bigl\{\sum_{z\in\Z^2} \ell_n(z)\xi(z)\geq b_n\Bigr\}\leq
\exp\Bigl\{-\frac{(1-\epsilon/2)b_n^2}{2\sigma^2 \ell_n^{\ssup 2}}\Bigr\}\leq
\exp\Bigl\{-\frac{(1-\epsilon)b_n^2}{2\sigma^2(i+1)\delta a_n}\Bigr\},
$$
and thus, applying (\ref{p2}) again and recalling the definition of $a_n$, for sufficiently large~$n$,
\begin{equation}\begin{aligned}\label{p6}
\mathbb{P}& \,\Bigl\{X_n\geq b_n, \,\ell_n^{\ssup \infty}\leq\eta i \delta \sqrt{n}, \, \ell_n^{\ssup 2}\in(i\delta a_n,(i+1)\delta a_n]\Bigr\}\\
&\leq\exp\Bigl\{-\frac{(1-\epsilon)b_n^2}{2\sigma^2(i+1)\delta a_n}\Bigr\} \,
\mathbb{P}\{\ell_n^{\ssup 2}>i\delta a_n\}
\leq\exp\Bigl\{-\frac{(1-\epsilon)b_n}{2\sigma^2(i+1)\delta \sqrt{n}}
-\frac{(1-\epsilon)\sqrt{\det\Gamma}i\delta b_n}{2\varkappa^4 \sqrt{n}}\Bigr\}.
\end{aligned}\end{equation}
It remains to consider the summand corresponding to $i=0$ in \eqref{p3}, which for any
$\eta>0$ is bounded~by
\begin{equation}\label{p7}
\mathbb{P}\Bigl\{X_n\geq b_n, \, \ell_n^{\ssup 2} \le \delta a_n, \,\ell_n^{\ssup \infty}\leq\eta \sqrt{n} \Bigr\}
+\mathbb{P}\Big\{  \ell_n^{\ssup \infty}>\eta \sqrt{n}\Big\}\, .\end{equation}
Applying Chebyshev's inequality on the event $\{ \ell_n^{\ssup 2} \le \delta a_n, \, \ell_n^{\ssup \infty}\leq\eta \sqrt{n}\}$
we get, for any $a>0$ and $\eta<\theta/a$,
$$\mathsf{P}\Big\{  \sum_{z\in\Z^2} \ell_n(z)\xi(z)\geq b_n\Big\}
\le \exp\Big\{-a\, \frac{b_n}{\sqrt{n}} + C \, \sum_{z\in\Z^2}\frac{a^2}{n}\, \ell_n^2(z) \Big\},$$
for a constant $C>0$ depending only on the distribution of the scenery
and the random walk. Using this estimate for
$a=1/(4C\delta)$ and $\eta<4C\delta\theta$ we get
\begin{equation}\label{p7b}
\mathbb{P}\bigl\{X_n\geq b_n, \, \ell_n^{\ssup 2} \le \delta a_n, \,\ell_n^{\ssup \infty}\leq\eta \sqrt{n} \bigr\}\leq
\exp\Big\{ -\frac{b_n}{\sqrt{n}}\, \frac{1}{8\delta C} \Big\}\, .
\end{equation}
Combining (\ref{p3}) -- (\ref{p7b}) gives us
\begin{equation}\begin{aligned}\label{p8}
\mathbb{P}\{X_n\geq b_n\}\leq\sum_{i=1}^{N-1} &
\exp\Bigl\{-\frac{(1-\epsilon)b_n}{2\sigma^2(i+1)\delta \sqrt{n}}-\frac{(1-\epsilon)i\delta b_n \sqrt{\det\Gamma}}{2\varkappa^4 \sqrt{n}}\Bigr\}\\
& +\exp\Bigl\{-\frac{b_n}{\sqrt{n}} \, \frac{1}{8\delta C}\Bigr\} +N\,\exp\Bigl\{-c\frac{\eta \delta \sqrt{n} }
{\log n}\Bigr\} +\exp\Bigl\{-\frac{N\delta b_n \sqrt{\det\Gamma}}{2\varkappa^{4} \sqrt{n}}\Bigr\}.
\end{aligned}\end{equation}
It is easily seen, that
$$\begin{aligned}
\lim_{n\to\infty}\frac{\sqrt{n}}{b_n}\log\sum_{i=1}^{N-1}  &
\exp\Bigl\{-\frac{(1-\epsilon)b_n}{2\sigma^2(i+1)\delta \sqrt{n}}-\frac{(1-\epsilon)i\delta b_n\sqrt{\det\Gamma}}{2\varkappa^4 \sqrt{n}}\Bigr\} \\
& =-(1-\epsilon)\min_{1\leq i\leq N-1}\Bigl(\frac{1}{2\sigma^2(i+1)\delta}+\frac{i\delta\sqrt{\det\Gamma}}{2\varkappa^4}\Bigr).
\end{aligned}$$
Furthermore, if we choose $\delta>0$ small and $N$ large, we get
$$ \min_{1\leq i\leq N-1}\Bigl(\frac{1}{2\sigma^2(i+1)\delta}+\frac{i\delta\sqrt{\det\Gamma}}{2\varkappa^4}\Bigr)\geq
(1-\epsilon)\,\min_{x>0}\,\Bigl(\frac{1}{2\sigma^2 x}+\frac{x\sqrt{\det\Gamma}}{2\varkappa^4}\Bigr)=
(1-\epsilon)\frac{(\det\Gamma)^{1/4}}{\sigma\varkappa^2}.
$$
Therefore, for all $n$ large enough,
\begin{equation}
\label{p9}
\sum_{i=1}^{N-1}
\exp\Bigl\{-\frac{(1-\epsilon)b_n}{4\sigma^2(i+1)\delta n^{1/2}}-\frac{(1-\epsilon)i
\delta b_n\sqrt{\det\Gamma} }{\varkappa^4  \sqrt{n}}\Bigr\}
\leq\exp\Bigl\{-(1-\epsilon)^3\frac{b_n\,(\det\Gamma)^{1/4}}{\sigma\varkappa^2 \sqrt{n}}\Bigr\}.
\end{equation}
Making first $\delta$ smaller, and then $N$ larger, if necessary, we see that all other terms
in \eqref{p8} are of smaller order than \eqref{p9}. Taking into account that
$\epsilon>0$ was arbitrary, we have
$$ \limsup_{n\to\infty}\frac{\sqrt{n}}{b_n}\log \prob\{X_n\geq b_n\}
\leq-\frac{(\det\Gamma)^{1/4}}{\sigma\varkappa^2}.$$

To obtain a \emph{lower bound}, note that for all $0<\mu<\lambda$ and $\eta>0$,
\begin{equation}
\label{p12}
\mathbb{P}\{X_n\geq b_n\}\geq\mathbb{P}\big\{X_n\geq b_n, \, \ell_n^{\ssup 2}\in[\mu a_n,\lambda a_n],
\, \ell_n^{\ssup \infty}\leq\eta \sqrt{n} \big\}\, ,
\end{equation}
where we still use $a_n=b_n\sqrt{n}$. Recall \eqref{4.1.'} and the definition of $L_n$
and $V_n$. Note that on the set $\{\ell_n^{\ssup 2}\in[\mu a_n,\lambda a_n], \,\ell_n^{\ssup \infty}\leq\eta \sqrt{n}\}$
and for sufficiently large~$n$, we  have $\frac 32 V_n \le \frac 32 \sigma(\lambda b_n)^{1/2} n^{1/4} \le b_n
\le \eta a_n / \ell_n^{\ssup \infty}
\le ( \eta/\mu\sigma^2) V_n^2/ \ell_n^{\ssup \infty} \le V_n/(196 L_n)$
if $\eta>0$ is sufficiently small. Hence,
$${\sf P}\Big\{ \sum_{z\in \Z^d} \ell_n(z) \xi(z) \ge b_n \Big\}
\ge \big(1-\Phi(\sfrac{b_n}{V_n})\big)\, \exp\big\{ -c_1 b_n^3\, L_n V_n^{-3} \big\}
\, \big( 1- c_2 b_n L_n V_n^{-1} \big) .$$
We observe that $L_n\le \frac{\gamma\eta}{\sigma^2}\, \frac{\sqrt{n}}{V_n}$ and hence
$$b_n^3 L_n V_n^{-3} \le \frac{\gamma\eta}{\sigma^2}\, b_n^3 \, \frac{\sqrt{n}}{V_n^4} \le
\frac{\gamma\eta}{\sigma^6\mu^2} \, \frac{b_n}{\sqrt{n}}
\qquad\mbox{ and } \qquad
b_nL_nV_n^{-1} \le \frac{\gamma\eta}{\sigma^2}\,  b_n \, \frac{\sqrt{n}}{V_n^2}
\le \frac{\gamma\eta}{\mu\sigma^4}\, .$$
Therefore, for all large~$n$,
\begin{eqnarray}
\label{p13}
\nonumber
\mathbb{P}\{X_n\geq b_n\}\geq\exp\Bigl\{-(1+\epsilon)\frac{b_n}{2\mu\sigma^2 \sqrt{n}}-\frac{c_1\gamma\eta b_n}{\mu^2\sigma^6\,\sqrt{n}}\Bigr\}
\Bigl(1-c_2\,\frac{\gamma\eta}{\mu\sigma^4}\Bigr)\\
\times\Bigl[\mathbb{P}\big\{\ell_n^{\ssup 2}\in[\mu a_n,\lambda a_n])\big\}-\mathbb{P}\big\{\ell_n^{\ssup \infty}>\eta\sqrt{n} \big\}\Bigr].
\end{eqnarray}
{F}rom (\ref{p2}) we conclude that for all $\mu<\lambda$,
\begin{equation}
\label{p14}
\log\mathbb{P}\big\{\ell_n^{\ssup 2}\in[\mu a_n,\lambda a_n]\big\}\sim-\frac{\mu b_n \sqrt{\det\Gamma}} {2\varkappa^4n^{1/2}}.
\end{equation}
Applying (\ref{p14}) and (\ref{p5}) to the right hand side of (\ref{p13}), we get for
$n^{1/2}\log n\ll b_n\ll n/\log n$,
$$\liminf_{n\to\infty}\frac{\sqrt{n}}{b_n}\log\mathbb{P}\{X_n\geq b_n\}\geq
-\frac{1+\epsilon}{2\mu\sigma^2}-\frac{c_1\gamma\eta}{\mu^2\sigma^6}-\frac{\mu\sqrt{\det\Gamma}}{2\varkappa^4}.$$
Since $\epsilon,\eta>0$ can be chosen arbitrarily small, and $\mu$ is arbitrary,
$$\liminf_{n\to\infty}\frac{\sqrt{n}}{b_n}\log\mathbb{P}\{X_n\geq b_n\}\geq
-\min_{\mu>0}\Bigl(\frac{1}{2\mu\sigma^2}+\frac{\mu\sqrt{\det\Gamma}}{2\varkappa^4}\Bigr)
=-\frac{(\det\Gamma)^{1/4}}{\sigma\varkappa^2}.$$ This completes the
proof of Theorem~\ref{maind=2}(b).\qed

\subsection{Proof of Theorem~\ref{maind=2}(c)}

We now assume that $b_n:=a\sqrt{n}\log n$.
In this case we use the following decomposition,
\begin{eqnarray}
\label{p15}
\nonumber
\mathbb{P}\{X_n\geq b_n\} & \leq &
\mathbb{P}\big\{X_n\geq b_n, \,\ell_n^{\ssup \infty}\leq\gamma_n,\, \ell_n^{\ssup 2}-\mathbb{E}\ell_n^{\ssup 2}
\leq\delta a_n\big\}\\ \nonumber & &
+\,\sum_{i=1}^N\mathbb{P}\big\{X_n\geq b_n, \,\ell_n^{\ssup \infty}\leq\gamma_n,\,\ell_n^{\ssup 2}-
\mathbb{E}\ell_n^{\ssup 2}\in(i\delta a_n,(i+1)\delta a_n]\big\}\\
& & +\,\mathbb{P}\big\{\ell_n^{\ssup \infty}>\gamma_n\big\}+\mathbb{P}\big\{\ell_n^{\ssup 2}
-\mathbb{E}\ell_n^{\ssup 2}>N\delta a_n\big\},\nonumber
\end{eqnarray}
here $a_n:=n\log n$, $\gamma _n:=\eta n\log n/b_n$. Estimating every
term as in the proof of the upper bound in (b) and using the
relation $\mathbb{E}\ell_n^{\ssup
2}\sim(\pi\sqrt{\det\Gamma})^{-1}n\log n$, one can get
$$
\limsup_{n\to\infty}\frac{1}{\log n}\mathbb{P}\{X_n\geq b_n\}\leq
-\min_{x\geq0}\Bigl(\frac{a^2}{2\sigma^2((\pi\sqrt{\det\Gamma})^{-1}+x)}-\frac{x\sqrt{\det\Gamma}}{2\varkappa^4}\Bigr)=I(a).
$$
In order to get a lower bound we consider the cases $a\leq\sigma/(\pi\varkappa^2(\det\Gamma)^{1/4})$ and
$a>\sigma/(\pi\varkappa^2(\det\Gamma)^{1/4})$ separately. In the first case we use
$$
\mathbb{P}\{X_n\geq b_n\}\geq \mathbb{P}\Big\{X_n\geq b_n,\, \ell_n^{\ssup \infty}\leq\gamma_n,\,
|\ell_n^{\ssup 2}-\mathbb{E}\ell_n^{\ssup 2}|\leq\delta a_n\Big\}, $$
and in the second case
$$
\mathbb{P}\{X_n\geq b_n\}\geq \mathbb{P}\Big\{X_n\geq b_n,\, \ell_n^{\ssup \infty}\leq\gamma_n, \,
\ell_n^{\ssup 2}-\mathbb{E}\ell_n^{\ssup 2}\in(\mu a_n,\lambda a_n]\Big\}
$$
for some $0<\mu<\lambda$. The further proof is similar to that of
the lower bound in Theorem~\ref{maind=2}(b) and details are left to
the reader.\qed

\section{Large deviations in dimension $d=2$: Proof of Proposition~\ref{special}}

We first derive an \emph{upper bound} for $\prob\{X_n\geq b_n\}$.
For arbitrary $N\geq1$ and $0<\delta<1$,
\begin{equation}
\label{P3}
\prob\{X_n\geq b_n\}\leq\sum_{i=0}^{N-1}\prob\{X_n\geq b_n,\ \ell_n^{\ssup \infty}\in(i\delta a_n,(i+1)\delta a_n]\}
+\prob\{\ell_n^{\ssup \infty}\geq \delta N a_n\},
\end{equation}
where $a_n:=(b_n\log n)^{1/2}$.
By assumption~(\ref{P1}), there exists $C_\delta$ such that
$$\me e^{h\xi(0)} \le \exp\{C_\delta h^2\}\quad\text{for}\quad h\leq(1-\delta)D.$$
{F}rom this bound and Chebyshev's inequality we get
\begin{equation}
\label{P4}
{\mathsf P}\Bigl\{\sum_{z\in\Z^d} \ell_n(z) \xi(z) \ge b_n\Bigr\}\leq
\exp\bigl\{-hb_n+C_\delta h^2 \ell_n^{\ssup 2} \bigr\}\quad\text{for}\quad h\leq(1-\delta)D/\ell_n^{\ssup \infty}.
\end{equation}
Letting here $h=(1-\delta)D/\ell_n^{\ssup \infty}$, we obtain
$$ {\mathsf P}\Bigl\{\sum_{z\in\Z^d} \ell_n(z) \xi(z) \ge b_n\Bigr\}\leq
\exp\Bigl\{-\sfrac{(1-\delta)D b_n}{\ell_n^{\ssup \infty}}\,
\bigl(1-\sfrac{C_\delta(1-\delta)D}{\ell_n^{\ssup \infty}b_n}\, \ell_n^{\ssup 2}\bigr)\Bigr\}.
$$
Therefore, for any $i\geq1$,
\begin{equation}\label{P5}\begin{aligned}
\prob\{X_n\geq b_n,& \ \ell_n^{\ssup \infty}\in(i\delta a_n,(i+1)\delta a_n]\} \\ & \leq
\exp\bigl\{-\sfrac{(1-\delta)^2D b_n}{(i+1)\delta a_n}\bigr\}\, \prob\{\ell_n^{\ssup \infty}>i\delta a_n\}+
\prob\bigl\{ \ell_n^{\ssup 2}>\sfrac{i\delta^2}{(1-\delta)C_\delta D}b_na_n\bigr\}.
\end{aligned}
\end{equation}
Using [GHK06, Lemma 1.3] and recalling the definition of $a_n$, we get
\begin{equation}
\label{P5'}
\log\prob\{\ell_n(0)>x a_n\}\sim -K_2x\frac{(b_n\log b_n)^{1/2}}{\log n-(1/2)\log b_n}
\sim-\frac{2K_2x}{2-\beta}\Bigl(\frac{b_n}{\log n}\Bigr)^{1/2}.
\end{equation}
Hence, arguing as in Lemma~\ref{maxdev}, for all $x\ge\delta$ and $n$ large enough $n$,
\begin{equation}
\label{P6}
\prob\{\ell_n^{\ssup \infty}>x a_n\}\leq\exp\Bigl\{
-(1-\delta)^2\frac{2K_2x}{2-\beta}\Bigl(\frac{b_n}{\log n}\Bigr)^{1/2}\Bigr\}.
\end{equation}
Combining (\ref{P5}) and (\ref{P6}), and noting that $b_n/a_n=(b_n/\log n)^{1/2}$, we obtain
\begin{equation}\label{P7}\begin{aligned}
\prob\{X_n\geq b_n,\ & \ell_n^{\ssup \infty}\in(i\delta a_n,(i+1)\delta a_n]\} \\
& \leq \exp\Bigl\{-\bigl(\sfrac{(1-\delta)^2D}{(i+1)\delta}+(1-\delta)^2\sfrac{2K_2i\delta}{2-\beta}\bigr)
\bigl(\sfrac{b_n}{\log n}\bigr)^{1/2}\Bigr\}+
\prob\bigl\{ \ell_n^{\ssup 2}>\sfrac{i\delta^2}{(1-\delta)C_\delta D}b_na_n\bigr\}.
\end{aligned}\end{equation}
Now we consider the probability corresponding to $i=0$. As $\ell_n^{\ssup \infty}\leq \delta a_n$, we can use
$h=(\delta^{-1}-1)Da_n^{-1}$ in (\ref{P4}). This gives us the bound
$$ {\mathsf P}\Bigl\{\sum_{z\in\Z^d} \ell_n(z) \xi(z) \ge b_n\Bigr\}\leq
\exp\Bigl\{-\sfrac{(1-\delta)D b_n}{\delta a_n}
\bigl(1-\sfrac{C_\delta(1-\delta)D}{\delta a_nb_n}\,\ell_n^{\ssup 2}\bigr)\Bigr\}. $$
Averaging over the random walk, we have
\begin{equation}
\label{P8}
\prob\{X_n\geq b_n,\ \ell_n^{\ssup \infty}\leq \delta a_n\}\leq
\exp\bigl\{-\sfrac{(1-\delta)^2D b_n}{\delta a_n}\bigr\}+
\prob\bigl\{ \ell_n^{\ssup 2} >\sfrac{\delta^2}{(1-\delta)C_\delta D}\,b_na_n\bigr\}.
\end{equation}
Applying (\ref{P6}) we obtain
\begin{equation}
\label{P9}
\prob\{\ell_n^{\ssup \infty}\geq \delta N a_n\}\leq\exp\Bigl\{
-c\delta N\bigl(\sfrac{b_n}{\log n}\bigr)^{1/2}\Bigr\}.
\end{equation}
Substituting (\ref{P7}) -- (\ref{P9}) into (\ref{P3}) gives
\begin{equation}\begin{aligned}
\label{P10}
\prob\{X_n\geq b_n\}& \leq\sum_{i=0}^{N-1}
\exp\Bigl\{-(1-\delta)^2\Bigl(\frac{}{(i+1)\delta}+\frac{2K_2i\delta}{2-\beta}\Bigr)\Bigl(\frac{b_n}{\log n}\Bigr)^{1/2}\Bigr\}\\
& +N\, \prob\bigl\{ \ell_n^{\ssup 2}>\sfrac{\delta^2}{(1-\delta)C_\delta D}\, b_na_n\bigr\}+
\exp\bigl\{-c\delta N\bigl(\sfrac{b_n}{\log n}\bigr)^{1/2}\bigr\}.
\end{aligned}\end{equation}
It is easily seen that
\begin{eqnarray*}
\lim_{n\to\infty}\Bigl(\frac{\log n}{b_n}\Bigr)^{1/2}\log \sum_{i=0}^{N-1}
\exp\Bigl\{-(1-\delta)^2\Bigl(\frac{D}{(i+1)\delta}+\frac{2K_2i\delta}{2-\beta}\Bigr)\Bigl(\frac{b_n}{\log n}\Bigr)^{1/2}\Bigr\}\\
=-(1-\delta)^2\min_{0\leq i<N}\Bigl(\frac{D}{(i+1)\delta}+\frac{2K_2i\delta}{2-\beta}\Bigr).
\end{eqnarray*}
Further, for small $\delta$ and large $N$ we have the inequality
$$\min_{0\leq i<N}\Bigl(\frac{D}{(i+1)\delta}+\frac{2K_2i\delta}{2-\beta}\Bigr)\geq
(1-\delta)\min_{x>0}\Bigl(\frac{D}{x}+\frac{2K_2x}{2-\beta}\Bigr)=
(1-\delta)\Big(\frac{8\,K_2 D}{2-\beta} \Big)^{1/2}.$$
Consequently, for all $n$ large enough,
\begin{equation}
\label{P11}
\sum_{i=0}^{N-1}
\exp\Bigl\{-(1-\delta)^2\Bigl(\frac{}{(i+1)\delta}+\frac{2K_2i\delta}{2-\beta}\Bigr)\Bigl(\frac{b_n}{\log n}\Bigr)^{1/2}\Bigr\}
\leq\exp\Bigl\{-(1-\delta)^4\Big(\frac{8\,K_2 D}{2-\beta} \Big)^{1/2}\Bigl(\frac{b_n}{\log n}\Bigr)^{1/2}\Bigr\}.
\end{equation}
Making $N$ larger, we see that the last term in (\ref{P10}) is of smaller order than (\ref{P11}).
By \eqref{p1} we obtain, for some constant~$c>0$,
$$
\log\prob\bigl\{\ell_n^{\ssup 2}>tb_na_n\bigr\}\sim -ct\,\Bigl(\frac{a_nb_n}{n}\Bigr).
$$
By our assumption, $b_n\log n\gg n$. Therefore, $n^{-1}a_nb_n=n^{-1}b_n^{3/2}\log^{1/2}n\gg (b_n/\log n)^{1/2}$.
This means that the probability term in (\ref{P10}) is negligible compared to (\ref{P11}). As a result we have
\begin{equation}
\label{P11'}
\limsup_{n\to\infty}\Bigl(\frac{\log n}{b_n}\Bigr)^{1/2}\log\prob\{X_n\geq b_n\}
\leq-(1-\delta)^4\Big(\frac{8\,K_2 D}{2-\beta} \Big)^{1/2}.
\end{equation}

To derive a \emph{lower bound} we note that
$${\mathsf P}\Bigl\{\sum_{z\in\Z^d} \ell_n(z) \xi(z) \ge b_n\Bigr\}\geq
{\mathsf P}\bigl\{\ell_n(0)\xi(0)\geq(1+\delta)b_n\bigl\}\,
{\mathsf P}\Bigl\{\sum_{z\neq 0} \ell_n(z) \xi(z) \ge -\delta b_n\Bigr\}.$$
Applying Chebyshev's inequality with second moments gives us
$$ {\mathsf P}\Bigl\{\sum_{z\in\Z^d} \ell_n(z) \xi(z) \ge b_n\Bigr\}\geq
{\mathsf P}\Bigl\{\ell_n(0)\xi(0)\geq(1+\delta)b_n\Bigl\} \,
\Bigl(1-\frac{\sigma^2\ell_n^2}{\delta^2 b_n^2}\Bigr). $$
Consequently,
\begin{equation}
\label{P12}
\prob\{X_n\geq b_n\}\geq (1-\delta)\, \prob\{\ell_n(0)\xi(0)\geq (1+\delta)b_n\}-
\prob\{\ell_n^{\ssup 2}>\delta^2\,b_n^2/\sigma^2\}.
\end{equation}
{F}rom (\ref{P1}) and (\ref{P5'}) we get, for every $x>0$,
$$\begin{aligned}
\prob\{\ell_n(0)\xi(0)\geq (1+\delta)b_n\}\geq\prob\{\ell_n(0)>xa_n\}
\geq \exp\Bigl\{-(1+\delta)^2\Bigl(\frac{D}{x}+\frac{2K_2x}{2-\beta}\Bigr)
\Bigl(\frac{b_n}{\log n}\Bigr)^{1/2}\Bigr\}.
\end{aligned}$$
Minimizing over $x$, we see that
\begin{equation}
\label{P13}
\prob\{\ell_n(0)\xi(0)\geq (1+\delta)b_n\}\geq
\exp\Bigl\{-(1+\delta)^2\Big(\frac{8\,K_2 D}{2-\beta} \Big)^{1/2}\Bigl(\frac{b_n}{\log n}\Bigr)^{1/2}\Bigr\}.
\end{equation}
As in the proof of the upper bound one can show that the last term in (\ref{P12}) is of smaller order than the
right hand side in (\ref{P13}). Therefore,
\begin{equation}
\label{P14}
\liminf_{n\to\infty}\Bigl(\frac{\log n}{b_n}\Bigr)^{1/2}\log\prob\{X_n\geq b_n\}
\geq-(1+\delta)^2\Big(\frac{8\,K_2 D}{2-\beta} \Big)^{1/2}.
\end{equation}
Combining (\ref{P11'}) and (\ref{P14}), and taking into account that $\delta$ is arbitrary, we get (\ref{P2}).
\qed

\vspace{0.5cm}

{\bf Acknowledgements:} We thank Amine Asselah and Fabienne Castell for interesting discussions and
for drawing our attention to the Preprint~\cite{As06}. This work is supported by an Advanced Research Fellowship
of the second author and by grants {f}rom DFG~(Germany) and EPSRC~(United Kingdom).

\vspace{1cm}


 {\footnotesize
\smallskip
\tableofcontents}

\end{document}